\documentclass{amsart}

\usepackage{latexsym,enumerate}
\usepackage{amsmath,amsthm,amsopn,amstext,amscd,amsfonts,amssymb}
\usepackage[ansinew]{inputenc}
\usepackage{verbatim}

\newcommand{\CC}{\mathbb{C}}
\newcommand{\NN}{\mathbb{N}}
\newcommand{\PP}{\mathbb{P}}
\newcommand{\RR}{\mathbb{R}}
\newcommand{\ZZ}{\mathbb{Z}}

\newcommand{\C}{\mathcal{C}}
\newcommand{\M}{\mathcal{M}}

\newcommand{\K}{\mathcal{K}}

\newtheorem{definition}{\sc Definition}[section]
\newtheorem{teo}{\sc Theorem}[section]
\newtheorem{prop}{\sc Proposition}[section]
\newtheorem{lemma}{\sc Lemma }[section]

\renewcommand{\span}[1]{\mbox{{\rm span\/}}#1}
\newcommand{\supp}[1]{\mbox{{\rm supp\/}}#1}

\renewcommand{\a}{\alpha}
\renewcommand{\b}{\beta}

\renewcommand{\d}{\delta}
\newcommand{\g}{\gamma}
\newcommand{\G}{\Gamma}
\renewcommand{\l}{\lambda}
\renewcommand{\L}{\Lambda}
\renewcommand{\O}{\Omega}

\newcommand{\s}{\sigma}
\newcommand{\z}{\zeta}

\begin{document}


\title[Sobolev spaces and
zeros of Sobolev orthogonal polynomials]{Sobolev spaces with respect to measures in
curves and zeros of Sobolev orthogonal polynomials}

\author{JOSE M. RODRIGUEZ$^{(1)}$, JOSE M. SIGARRETA$^{(2)}$}

\date{January 30, 2008.\\
$(1)\,\,\,$
Supported in part by three grants from M.E.C.\
(MTM 2006-13000-C03-02, MTM 2006-11976 and MTM 2006-26627-E), Spain,
and by a grant from
U.C.III$\,$M./C.A.M.\ (CCG06-UC3M/EST-0690), Spain.  \\
$(2)\,\,\,$
Supported in part by a grant from M.E.C.\
(MTM 2006-13000-C03-02), Spain,
and by a grant from
U.C.III$\,$M./C.A.M.\ (CCG06-UC3M/EST-0690), Spain.
}

\begin{abstract}
In this paper we obtain some practical criteria to bound the
multiplication operator in Sobolev spaces
with respect to measures in curves.
As a consequence of these results, we characterize the weighted Sobolev
spaces with bounded multiplication operator,
for a large class of weights. To have
bounded multiplication operator has important consequences in
Approximation Theory: it implies the uniform bound of the zeros of
the corresponding Sobolev orthogonal polynomials, and this fact
allows to obtain the asymptotic behavior of Sobolev orthogonal
polynomials. We also obtain some non-trivial results about
these Sobolev spaces with respect to measures;
in particular, we prove a main result in the theory: they are Banach spaces.
\end{abstract}

\maketitle{}

\

\centerline{Running headline: Sobolev spaces and zeros of Sobolev orthogonal
polynomials}

\

\centerline{1991 AMS Subject Classification: 41A10, 46E35, 46G10}

\

\centerline{Key words: Multiplication operator; location of zeros}

\centerline{Sobolev orthogonal polynomials; weight; weighted Sobolev spaces on curves}

\

\

Jos\'e M. Rodr{\'\i}guez

Departamento de Matem\'aticas

Escuela Polit\'ecnica Superior

Universidad Carlos III de Madrid

Avenida de la Universidad, 30

28911 Legan\'es (Madrid), SPAIN

Telephone number: 34 91 6249098

Fax number: 34 91 6249151

email: {\tt jomaro@math.uc3m.es}

\

Jos\'e M. Sigarreta

Departamento de Matem\'aticas

Universidad Popular Aut\'onoma del Estado de Puebla

21 Sur 1103, Col. Santiago

CP 72410, Puebla, MEXICO.

email: {\tt jsigarre@math.uc3m.es}

\newpage

\section{Introduction.}

Weighted Sobolev spaces are an interesting topic in many
fields of Mathematics. In the classical book \cite{Ku}
we can find the point of view of Partial Differential Equations.
We are interested in
the relationship between this topic and Approximation Theory in
general, and Sobolev orthogonal polynomials in particular.

Sobolev orthogonal polynomials have been more and more investigated in recent years.
In particular, in \cite{IKNS1} and \cite{IKNS2}, the authors showed that
the expansions with Sobolev orthogonal polynomials can avoid the Gibbs
phenomenon which appears with classical orthogonal series in $L^2$.
The papers \cite{APRR}, \cite{BFM}, \cite{CM}, \cite{FMP}, \cite{LPP} and
\cite{RY} deal with Sobolev spaces on curves and more
general subsets of the complex plane.

Sobolev orthogonal polynomials on the unit circle and, more generally, on curves is a
topic of recent and increasing interest in approximation theory; see, for instance,
\cite{CM} and \cite{FMP} (for the unit circle) and \cite{M-F} and \cite{BFM} (for the
case of Jordan curves).

In \cite{RARP1}, \cite{RARP2}, \cite{R1}, \cite{R2} and \cite{R3}
the authors solved the following specific problems:

1) Find hypotheses on
general measures $\mu=(\mu_0, \mu_1, \dots, \mu_k)$ in $\RR$,
as general as possible, so that
we can define a Sobolev space $W^{k,p}(\mu)$ whose elements are functions.
These measures are called $p$-admissible.

2) If a Sobolev norm with general measures
$\mu=(\mu_0, \mu_1, \dots, \mu_k)$ in $\RR$ is
finite for any polynomial, what is the completion, $\PP^{k,p}(\mu)$,
of the space of polynomials with respect to the norm in $W^{k,p}(\mu)$?
This problem has been studied previously in some
particular cases (see e.g. \cite{ELW1}, \cite{EL}, \cite{ELW2}).

In \cite{APRR} and \cite{RY} these results are extended to weighted
Sobolev spaces on curves in the complex plane.

One of the central problems in the theory of Sobolev orthogonal
polynomials is to determine its asymptotic behavior. In \cite{LP}
the authors show how to obtain the $n$-th root asymptotic of Sobolev
orthogonal polynomials if the zeros of these polynomials are
contained in a compact set of the complex plane. Although the
uniform bound of the zeros of orthogonal polynomials holds for every
measure with compact support in the case without derivatives
($k=0$), it is an open problem to bound the zeros of Sobolev
orthogonal polynomials. The boundedness of the zeros is a
consequence of the boundedness of the multiplication operator $\M
f(z)=z\,f(z)$ in the corresponding space $\PP^{k,2}(\mu)$: in fact,
the zeros of the Sobolev orthogonal polynomials are contained in the
disk $\{z:\, |z|\le \|\M\|\}$ (see \cite{LPP}).

In \cite{RARP2}, \cite{R2}, \cite{APRR} and \cite{R4},
there are some answers to the
question stated in \cite{LP} about some
conditions for $\M$ to be bounded.

The main aim of this paper is to find conditions
(which should be easy to check in practical cases)
implying the boundedness of these zeros,
when the supports of the measures are contained
in a curve in the complex plane
(see Theorems \ref{mult1}, \ref{mult2}, \ref{mult3} and \ref{mult4}).
Theorem \ref{mult1} is a general result which can be applied for
a wide class of weights: in fact,
Theorems \ref{mult2}, \ref{mult3} and \ref{mult4} are consequences of it.
In particular, Theorem \ref{mult3} states the following characterization:
If $d\mu_j=w_j ds$ and $w_j$ is piecewise monotone for $1\le j \le k$,
then $\M$ is bounded if and only if $\K(\mu)=0$
(this condition means that the Sobolev norm
$\|f\|_{W^{k,p}(\mu)}:= (\sum_{j=0}^k
\|f^{(j)}\|^p_ {L^p(\mu_j)} )^{1/p}$
``is a norm"; see Definition \ref{d:K} for the precise definition of $\K(\mu)$).
Condition $\K(\mu)=0$ is easy to check in practical cases
(see Propsition \ref{k=1} for a characterization if $k=1$, and Theorem \ref{mult4} and
the Remark after Theorem \ref{mult1} for some sufficient conditions for any $k$).
The hypothesis about the monotonicity of $w_j$ is a weak one, since it is verified
in almost every example (for instance, every Jacobi, Jacobi-Angelesco and Polacheck
weight satisfies it). Theorem \ref{mult2} is a generalization of Theorem \ref{mult3}
with $d\mu_j=w_j ds+ d(\mu_j)_s$.
Theorem \ref{mult4} deals with weights ``similar" to some power, in some sense.

These results are new for Sobolev orthogonal polynomials in curves,
and even for Sobolev orthogonal polynomials in the real line.
Theorem \ref{mult3} is an improvement of \cite[Theorem 4.3]{R4}:
in \cite{R4} appears a different Sobolev space (in an interval
$I$), denoted by $W_{ko}^{k,p}(I,\mu)$, verifying
$\PP^{k,p}(I,\mu) \subseteq W^{k,p}(I,\mu) \subseteq W_{ko}^{k,p}(I,\mu)$.
Since we usually have
$\PP^{k,p}(I,\mu) = W^{k,p}(I,\mu) \neq W_{ko}^{k,p}(I,\mu)$
(see \cite[Theorem 6.1]{APRR}),
it is obvious that it is better to work with $W^{k,p}(I,\mu)$ in order
to obtain results about the multiplication operator in $\PP^{k,p}(I,\mu)$.
Furthermore, we can define $W^{k,p}(I,\mu)$ for any $p$-admissible measure $\mu$
(see Definition \ref{8}), but we can define $W_{ko}^{k,p}(I,\mu)$ just for a smaller
set of measures (see \cite[Definition 2.3]{R4}).

(The advantage of $W_{ko}^{k,p}(I,\mu)$ is that it can be defined in a simpler and
faster way than $W^{k,p}(I,\mu)$.)

We have proved also some technical results
about weighted Sobolev spaces.
The main technical result is
Theorem \ref{t:c1}, which says that $W^{k,p}(\mu)$ is a Banach space for every
$p$-admissible measure $\mu$.
This central result in the theory of Sobolev spaces has an interesting consequence
for the study of the multiplication operator:
if $\M$ is bounded in $W^{k,p}(\mu)$
and $\PP \subseteq W^{k,p}(\mu)$,
then it is bounded in $\PP^{k,p}(\mu)$,
since $W^{k,p}(\mu)$ is a complete space.
This is a crucial fact, since the results in this paper deal with
the multiplication operator in $W^{k,p}(\mu)$;
we need to work with $W^{k,p}(\mu)$ instead of $\PP^{k,p}(\mu)$,
since in $W^{k,p}(\mu)$ we have powerful tools, like
Theorems \ref{t:4.1} and \ref{t:8.3});
furthermore, the elements of $W^{k,p}(\mu)$ are functions
and the elements of $\PP^{k,p}(\mu)$ are equivalence classes
of sequences of plynomials.

The outline of the paper is as follows.
Sections 2 and 3 are dedicated to the definitions
and previous results which will be useful.
In Section 4 we obtain some improvements of the results in Section 3,
which are interesting by themselves and
simplify many results about weighted Sobolev spaces in
\cite{RARP1}, \cite{RARP2}, \cite{R1}, \cite{R2}, \cite{R3}, \cite{APRR} and \cite{RY}.
We prove the results on multiplication operator in Section 5.

In order to make easy the reading of the paper to those people mainly interested
in the boundedness of zeros of Sobolev orthogonal polynomials,
Section $5$ is almost self-contained,
and does not depend too much on the rest of the paper.
In this section, whenever a previous technical result is used,
there is a precise reference to its location.

\smallskip

Now we introduce the notation we use.

\smallskip

{\bf Notation.} If $A$ is a Borel set in a curve, $\chi_{_{A}}, \,\,
\sharp A$ and $\overline{A}$ denote, respectively,
the characteristic function, the cardinal and the closure
of $A$. By $f^{(j)}$ we mean the $j$-th distributional derivative of
$f$. $\PP$ denotes the set of polynomials and $\PP_n$ the set of
polynomials of degree less or equal than $n$. We say that an
$n$-dimensional vector satisfies a one-dimensional property if each
coordinate satisfies this property. Finally, the constants in the
formulae can vary from line to line and even in the same line.

\smallskip

{\bf Acknowledgements.}
I would like to thank Professor Dmitry V. Yakubovich
for many useful discussions and several ideas which have improved
the presentation of the paper.

\

\section{Curves and derivatives along curves.}

In this section we introduce a definition of {\it derivative along a curve}
in the complex plane,
as an extension of the usual complex derivative,
which will need in the rest of the paper.
A detailed study of this concept (with the proofs of the results stated here)
can be found in \cite[Chapter 2]{APRR}.
Every curve will be simple, rectifiable and oriented.
Any closed curve is positively oriented (counter clockwise).

\begin{definition}
{\rm (a)}
If $\g$ is not a closed curve, its orientation give a natural order in $\g$;
then, if $z_0,z_1 \in \g$ and $z_0<z_1$, we can consider the \emph{arcs}
$[z_0,z_1]$, $(z_0,z_1)$, $[z_0,z_1)$ and $(z_0,z_1]$.
If $\g$ is a closed curve we denote by $[z_0,z_1]$ the arc of $\g$ joining $z_0$ with
$z_1$ in the positive sense; then we also have a natural order in each arc $[z_0,z_1]$.
We denote by $\int_{z_0}^{z_1} g(\zeta) d\zeta$ the
complex integral of the function $g$ along the arc $[z_0,z_1]$.

\smallskip

{\rm (b)} Let $z_0$ be a fixed point in $\g$. If $\g$ is compact we say that $f\in
AC^k(\g)$ if $f$ can be written as
\begin{equation}
\label{eq:0}
f(z)=q(z)+ \int_{z_0}^z h(\zeta) \, \frac {(z-\zeta)^{k-1}}{(k-1)!} \, d\zeta \,,
\end{equation}
for some $h\in L^1(\g,ds)$ and some polynomial $q\in \PP_{k-1}$.
If $\g$ is a closed curve we require also the function $h\in L^1(\g,ds)$
to verify $\int_\g h(\zeta)\,\z^i\,d\z=0$, for $0\le i<k$.
When $\g$ is not compact, we say that $f\in AC_{loc}^k(\g)$ if it
can be split as in $(\ref{eq:0})$ with $h\in L_{loc}^1(\g,ds)$.

\smallskip

{\rm (c)} If $f\in AC^k_{loc}(\g)$ and $z_0 \in \g$, we define its
derivative
$f'$ along $\g$ as
$$
f'(z)=q'(z) + \int_{z_0}^z h(\zeta) \, \frac {(z-\zeta)^{k-2}}{(k-2)!} \, d\zeta \,,
$$
where $q'(z)$ means the classical derivative of $q(z)$ and $\int_{z_0}^z
h(\zeta)(z-\zeta)^{-1}/(-1)! \,d\zeta$ means $h(z)$.
\end{definition}

Obviously, if $\g$ is a compact real interval, the space $AC^1(\g)$ is the set
of absolutely continuous functions in $\g$.
If $\g$ is a closed curve and $f\in AC^k(\g)$, we have
$\int_\g h(\zeta)\,(z-\z)^{k-1}\,d\z=0$
for every $z\in \g$. This property is equivalent to $f^{(j)}$ being continuous in
$\g$ for $0\le j<k$, where $f^{(j)}$ denotes the $j$-th derivative
(according to the previous definition) of $f$.

\smallskip

We also note that it is natural to define the derivative along $\g$ in this
way, since this is the ``inverse" of integration:
$$
\begin{aligned}
\int_{z_0}^z \int_{z_0}^\xi h(\zeta)\, \frac {(\xi -\zeta)^{k-2}}{(k-2)!} \, d\zeta \,
d\xi
& = \int_{z_0}^z \int_\zeta^z h(\zeta)\, \frac {(\xi -\zeta)^{k-2}}{(k-2)!} \, d\xi \,
d\zeta
\\
& = \int_{z_0}^z h(\zeta)\, \Big[ \frac {(\xi-\zeta)^{k-1}}{(k-1)!}
\Big]_{\xi=\zeta}^{\xi=z} \, d\z
= \int_{z_0}^z h(\zeta) \, \frac {(z-\zeta)^{k-1}}{(k-1)!} \, d\zeta \,.
\end{aligned}
$$

\noindent {\bf Remarks.}

{\bf 1.}
Note that if $f$ is holomorphic in a region containing
$\g$, then $f'$ is the usual complex derivative of $f$ at almost
every point of $\g$.

{\bf 2.}
If $f\in AC^k_{loc}(\g)$ and $f^{(j)}=0$ a.e. in $\g$, for some $0< j \le k$,
then $f\in \PP_{j-1}$.

\smallskip

It can be shown that this definition of derivative is independent of the
representation of $f$ we are using,
and that it verifies the properties of usual derivation:
linearity, Leibniz' rule, approximation by Taylor polynomials,...
(see \cite{APRR}).

\

\section{Background and previous results on Sobolev spaces.}

The main concepts that we need in order to state our
results are contained in the following definitions. The first one
is a class of weights that will be the absolutely continuous part
of our measures.

\begin{definition}
\label{2}
Given $1\le p < \infty$, a curve $\g$ and
a set $A$ which is a union of arcs in $\g$,
we say that a weight $w$ in $\g$
belongs to $B_p(A)$ if $w^{-1}\in L_{loc}^{1/(p-1)} (A)$
(if $p=1$, then $1/(p-1)=\infty$).
\end{definition}

It is possible to construct a similar theory with $p=\infty$.
We refer to
\cite{APRR}, \cite{PQRT1}, \cite{PQRT2} and \cite{PQRT3}
for the case $p=\infty$.

\smallskip

If the curve $\g$ is $\RR$, then $B_p(\RR)$ contains,
as a very particular case, the classical $A_p(\RR)$
weights appearing in Harmonic Analysis.
The classes $B_p(\O)$, with $\O\subseteq\RR^n$,
have been used in other definitions of
weighted Sobolev spaces on $\RR^n$ in \cite{KO}.

\smallskip

We  consider vectorial measures $\mu =(\mu_0, \dots , \mu_k)$ in
the definition of our Sobolev space in a curve. We assume that we
can make for each scalar measure the decomposition $d\mu_j=
d(\mu_j)_s + w_j ds$, where $(\mu_j)_s$ is singular with respect
to the Euclidean arc-length and $w_j$ is a non-negative Lebesgue
measurable function in $\g$ (by Radon-Nikodym's Theorem, we can
make this decomposition, for instance, if $\mu_j$ is $\s$-finite).

In \cite{KO}, Kufner and Opic define the following sets:

\begin{definition}
\label{4}
Let us consider $1\le p < \infty$ and a vectorial measure $\mu =
(\mu_0, \dots , \mu_k)$ in $\g$. For $0\le j \le k$ we define the open set
$$
\O_j:=\big\{ z \in \g\; : \; \exists \hbox{  an open neighbourhood  }
V \hbox{  of  } z  \hbox{ in the curve $\g$ with  } w_j\in B_p(V)\big\}\,.
$$
\end{definition}

Note that we always have $w_j \in B_p(\O_j)$ for any $0\le
j \le k$. In fact, $\O_j$ is the largest open set $U$ with $w_j\in B_p(U)$.
It is easy to check that  if $f^{(j)} \in L^p(\O_j,w_j)$ with $1\le j
\le k $, then $f^{(j)} \in L^1_{loc}(\O_j)$, and therefore
$f^{(j-1)} \in AC_{loc}^1(\O_j)$.

Since the precise definition of Sobolev space requires some technical concepts
(see Definition \ref{9}), we would like to introduce here a heuristic definition
of Sobolev space and an example which will help us to understand the technical process
that we will follow in order to reach to Definition \ref{9}.

\begin{definition}
\label{heuristic}
(Heuristic definition.)
Let us consider $1\le p< \infty$
and a $p$-admissible vectorial measure $\mu=(\mu_0,\dots,\mu_k)$
in $\g$. We define the \emph{Sobolev space}
$W^{k,p}(\g,\mu)$ as the space of equivalence classes of
$$
\begin{aligned}
V^{k,p}(\g,\mu):=\Big\{f:
&
\g\rightarrow\CC \; \;  /  \;\;
\big\|f\big\|_{W^{k,p}(\g,\mu)}:=\Big(\sum_{j=0}^k
\big\|f^{(j)}\big\|^p_ {L^p(\g,\mu_j)} \Big)^{1/p} < \infty \,,
\\
&
f^{(j)}\in AC_{loc}^1 (\O_{j+1}\cup\cdots\cup\O_k) \hbox{ and }
\, f^{(j)} \hbox{ satisfies ``pasting conditions" for }
0\le j < k \Big\}\,,
\end{aligned}
$$
with respect to the seminorm $\|\cdot \|_{W^{k,p}(\g,\mu)}$.
\end{definition}

\smallskip

These pasting conditions are natural: a function must be as
regular as possible. In a first step, we check if the functions and
its derivatives are absolutely continuous up to the boundary (this
fact holds in the following example), and then we join the contiguous
intervals:

\medskip

\noindent
{\bf Example.}
$\mu_0:=\d_0$, $\mu_1:=0$, $d\mu_2:=\chi_{{}_{\scriptstyle [-1,0]}}(x) dx$
and $d\mu_3:=\chi_{{}_{\scriptstyle [0,1]}}(x) dx$.

Since $\O_1=\emptyset$, $\O_2=(-1,0)$ and $\O_3=(0,1)$,
$W^{3,p}(\mu)$ is the space of equivalence classes of
$$
\begin{aligned}
V^{3,p}(\mu)
= \Big\{ f \ /  \
\|f\|_{W^{3,p}(\mu)}
&
< \infty \,,
\, \text{ $f$ satisfies ``pasting conditions",}
\\
&
\ \, f,f',\in AC((-1,0))
\, \text{ and }
\, f,f',f''\in AC((0,1))
\Big\}
\\
= \Big\{ f \ /  \
\|f\|_{W^{3,p}(\mu)}
&
< \infty \,,
\, \text{ $f$ satisfies ``pasting conditions",}
\\
&
\ \, f,f',\in AC([-1,0])
\, \text{ and }
\, f,f',f''\in AC([0,1])
\Big\}
\\
= \Big\{ f \ /  \
\|f\|_{W^{3,p}(\mu)}
&
< \infty \,,
\ \, f,f'\in AC([-1,1])
\, \text{ and }
\, f''\in AC([0,1])
\Big\}\,.
\end{aligned}
$$
In the current case, since $f$ and $f'$ are absolutely continuous in
$[-1,0]$ and in $[0,1]$, we require that both are absolutely
continuous in $[-1,1]$.

\smallskip

These heuristic concepts can be formalized as follows:

\begin{definition}
Let us consider $1\le p<\infty$, $\mu,\nu$ measures in
$\g$ and $z_0,z_1 \in \g$. We define
$$
\begin{aligned}
\L_{p,[z_0,z_1]}^+ (\mu, \nu) & := \sup_{z_0<z<z_1} \mu((z_0,z])
 \big\|(d\nu/ds)^{-1}\big\|_{L^{1/(p-1)}([z,z_1], ds)}\,,
\\
\L_{p,[z_0,z_1]}^- (\mu, \nu) & := \sup_{z_0<z<z_1} \mu([z,z_1))
 \big\|(d\nu/ds)^{-1}\big\|_{L^{1/(p-1)}([z_0,z], ds)}\,,
\end{aligned}
$$
where we use the convention $0\cdot \infty=0$.
\end{definition}

Before we state our theorems, let us state a version on curves of
a classical result in $\RR$ (see \cite{Mu2}, \cite{M}).
It will be generalized in Theorem \ref{t:4.1} below.

\

\noindent {\bf  Muckenhoupt inequality.} (\cite[Theorem 3.1]{APRR})
{\it Let us consider $1\le p<\infty$, $\mu_0,\mu_1$ measures in
$\g$ and $z_0,z_1 \in \g$. Then:

$(1)$ There exists a positive
constant $c$ such that
$$
\Big\|\int_{z}^{z_1} g(\zeta)\,d\z\Big\|_{L^p((z_0,z_1],\mu_0)}\le
 c\,
\big\|g\big\|_{L^p((z_0,z_1],\mu_1)}
$$
for any measurable function $g$ in $[z_0,z_1]$, if and only if
$\L_{p,[z_0,z_1]}^+ (\mu_0, \mu_1)<\infty$.

$(2)$ There exists a positive
constant $c$ such that
$$
\Big\|\int_{z_0}^z g(\z)\,d\z\Big\|_{L^p([z_0,z_1),\mu_0)}\le
 c\,
\big\|g\big\|_{L^p([z_0,z_1),\mu_1)}
$$
for any measurable function $g$ in $[z_0,z_1]$, if and only if
$\L_{p,[z_0,z_1]}^- (\mu_0, \mu_1)<\infty$.
}

\begin{definition}
Let us consider $1\le p<\infty$ and a curve $\g$. A vectorial
measure $\overline \mu=(\overline \mu_0,\dots,\overline \mu_k)$ is
a right completion of a vectorial measure
$\mu=(\mu_0,\dots,\mu_k)$ with respect to $z_0\in\g$ in a right
neighborhood $[z_0,z_1]$, if $\overline \mu_k=\mu_k$ in $\g$,
$\overline \mu_j=\mu_j$ in the complement of $(z_0,z_1]$ and
$$
\overline{\mu}_j=
\mu_j+\tilde{\mu}_{j}\,,\qquad  \text{in } (z_0,z_1] \;\, \text{
for } 0\le j<k\,,
$$
where $\tilde \mu_{j}$ is any measure satisfying
$\tilde \mu_j((z_0,z_1])<\infty$ and
$\Lambda_{p, [z_0,z_1]}^+ (\tilde \mu_j,\overline \mu_{j+1})<\infty$.
\end{definition}

The Muckenhoupt inequality guarantees that if $f^{(j)}\in L^p(\mu_j)$
and $f^{(j+1)}\in L^p(\overline\mu_{j+1})$, then $f^{(j)}\in L^p(\overline\mu_j)$.
If we work with absolutely continuous measures,
we also say that a vectorial weight $\overline{w}$ is a completion of $\mu$
(or of $w$). See some examples of completions in \cite{RARP1} and \cite{APRR}.

\medskip

\noindent
{\bf Remark.}
We can define a left completion of
$\mu$ with respect to $z_0$ in a similar way.

\begin{definition}
\label{d:regular}
For $1\le p< \infty$ and a vectorial measure $\mu$ in $\g$, we say that a
point $z_0\in \g$ is right $j$-regular $($respectively, left $j$-regular$)$,
if there exist a right completion $\overline\mu$
$($respectively, left completion$)$ of $\mu$
in $[z_0,z_1]$ and $j<i\le k$ such that
$\overline{w}_{i}\in B_p ([z_0,z_1])$
$($respectively, $B_p([z_1,z_0]))$.
Also, we say that a point $z_0\in \g$ is
$j$-regular, if it is right and left $j$-regular.
\end{definition}

\begin{definition}
\label{d:comparable}
We say that two functions $u,v$ are
\emph{comparable}
(and we write $u \asymp v$)
on the set $F\subseteq \g$
if there are positive constants $c_1,c_2$ such that
$c_1 v(z)\le u(z)\le c_2 v(z)$ for almost every $z\in F$.
\end{definition}

\noindent
{\bf Remarks.}

{\bf 1.} A point $z_0\in \g$ is
right $j$-regular $($respectively, left $j$-regular$)$, if at least one of the
following properties is verified:

\smallskip

{\rm (a)} \ There exist a right (respectively, left) neighborhood
$[z_0,z_1]$ (respectively, $[z_1,z_0]$) and $j<i\le k$ such that
$w_{i} \in B_p ([z_0,z_1])$ $($respectively, $B_p([z_1,z_0]))$.
Here we have chosen $\tilde w_j=0$.

{\rm (b)} \ There exist a right (respectively, left) neighborhood
$[z_0,z_1]$ (respectively, $[z_1,z_0]$) and $j<i\le k$, $\a>0$,
$\d< (i-j)p-1$, such that $w_{i}(z) \ge \a \, |z-z_0|^\delta$, for
almost every $z\in [z_0,z_1]$ $($respectively, $[z_1,z_0])$ and we
have $|z-z_0| \asymp |\g^{-1} (z) -\g^{-1}(z_0)|$ in $[z_0,z_1]$
(respectively, $[z_1,z_0]$), where $\g^{-1}$ denotes the inverse
of some parametrization with $\g' \in L^{\infty}([z_0,z_1])$. See
Lemma 3.4 in \cite{RARP1}.

{\bf 2.}
If $z_0$ is right $j$-regular (respectively, left), then it is also right
$i$-regular (respectively, left) for each $0\le i \le j$.

\medskip

When we use this definition we think of a point $\{z\}$ as the union of two
half-points $\{z^+\}$ and $\{z^-\}$. With this convention,
each one of the following sets
$$
\begin{aligned}
(z_0,z_1) \cup (z_1,z_2) \cup \{z_1^+\}
& = (z_0,z_1) \cup [z_1^+,z_2) \ne (z_0,z_2) \,,
\\
(z_0,z_1) \cup (z_1,z_2) \cup \{z_1^-\}
& = (z_0,z_1^-] \cup (z_1,z_2) \ne (z_0,z_2) \,,
\end{aligned}
$$
has two connected components, and the set
$$
(z_0,z_1) \cup (z_1,z_2) \cup \{z_1^-\} \cup \{z_1^+\} = (z_0,z_1) \cup (z_1,z_2) \cup \{z_1\} = (z_0,z_2)
$$
is connected.

We /use this convention in order to study the sets of continuity of functions:
we want that  if $f\in C(A)$ and $f\in C(B)$,
where $A$ and $B$ are union of arcs, then $f\in C(A\cup B)$.
With the usual definition of continuity in an arc,
if $f\in C([z_0,z_1))\cap C([z_1,z_2])$ then we do not have $f\in C([z_0,z_2])$. Of
course, we have $f\in C([z_0,z_2])$ if and only if
$f\in C([z_0,z_1^-])\cap C([z_1^+,z_2])$,
where by definition, $C([z_1^+,z_2])=C([z_1,z_2])$
and $C([z_0,z_1^-])=C([z_0,z_1])$. This idea can be formalized with a suitable
topological space.

\smallskip

Let us introduce some more notation.
We denote by $\O^{(j)}$ the set of $j$-regular points or half-points,
i.e., $z\in\O^{(j)}$  if and only if $z$ is $j$-regular,
we say that $z^+ \in \Omega^{(j)}$  if and only if $z$ is right $j$-regular, and
we say that $z^- \in \Omega^{(j)}$  if and only if $z$  is left $j$-regular.
Obviously, $\O^{(k)}=\emptyset$ and  $\Omega_{j+1}\cup \cdots \cup
\Omega_k\subseteq \O^{(j)}$.
Note that $\O^{(j)}$ depends on $p$ (see Definition 3.5).

\smallskip

\noindent {\bf Remark.}
If $0\le j <k$ and $J$ is an arc in $\g$,
$J\subseteq \O^{(j)}$, then the set
$J\setminus (\O_{j+1} \cup\cdots\cup\O_k)$ is discrete
(see the Remark before Definition 7 in \cite{RARP1}).
Consequently,
$\O^{(j)} \subseteq \overline{\O_{j+1} \cup\cdots\cup\O_k}$
for every $0\le j <k$.

\smallskip

Intuitively, $\O^{(j)}$ is the set of ``good" points at the level $j$
for the vectorial weight $(w_0,\dots , w_k)$ in $\g$:
every function $f$ in the Sobolev space must verify that
$f^{(j)}$ is continuous in $\O^{(j)}$.

\smallskip

Let us present now the class of measures that we use and the
definition of Sobolev space.

\begin{definition}
\label{8}
We say that the vectorial measure $\mu=(\mu_0,\dots,\mu_k)$ in $\g$ is
\emph{$p$-admissible} if $\mu_j^*(\g\setminus \Omega^{(j)})=0$, for
$1\le j < k$, and $\mu_k^* \equiv 0$,
where $d\mu_j^* := d\mu_j - w_j \chi_{_{\Omega_{j}}} \! ds$
(then $d\mu_k= w_k \chi_{_{\Omega_{k}}} \! ds$).

We say that the vectorial measure $\mu=(\mu_0,\dots,\mu_k)$ in $\g$ is
\emph{strongly $p$-admissible} if it is $p$-admissible and
$\supp (\mu_j^*|_{\L}) \subseteq \L$,
for any connected component $\L$ of $\Omega^{(j)}$
and $1\le j < k$.
\end{definition}

\noindent {\bf Remarks.}

{\bf 1.} The hypothesis of $p$-admissibility is natural. It would
not be reasonable to consider Dirac's deltas in $\mu_j$ in the
points where $f^{(j)}$ is not continuous.

{\bf 2.} Note that there is not any restriction on $\mu_0$.

{\bf 3.} Every absolutely continuous measure
$w=(w_0, \dots, w_k)$ with $w_j=0$
a.e. in $\g\setminus \Omega_j$ for every $1\le j \le k$,
is $p$-admissible and even strongly $p$-admissible
(since then $\mu_j^*=0$).
It is difficult to find a weight $w$ which does not satisfy this condition.

{\bf 4.} $(\mu_j)_s \le \mu_j^*$, and the equality usually holds.

\begin{definition}
\label{9}
Let us consider $1\le p< \infty$
and a $p$-admissible vectorial measure $\mu=(\mu_0,\dots,\mu_k)$
in $\g$. We define the \emph{Sobolev space}
$W^{k,p}(\g,\mu)$ as the space of equivalence classes
of
\begin{eqnarray}
V^{k,p}(\g,\mu):=\Big\{f:\g\rightarrow\CC \
&/& \ f^{(j)}\in AC_{loc}^1 (\O^{(j)}) \hbox{ for }
0\le j < k\hbox{
and }  \nonumber \\
&\,& \
\big\|f\big\|_{W^{k,p}(\g,\mu)}:=\Big(\sum_{j=0}^k
\big\|f^{(j)}\big\|^p_ {L^p(\g,\mu_j)} \Big)^{1/p}
 <\infty \Big\} \,, \nonumber
\end{eqnarray}
with respect to the seminorm $\|\cdot \|_{W^{k,p}(\g,\mu)}$.
\end{definition}

Perhaps this definition of Sobolev space is very technical, but it
has interesting properties: in many cases,
$\PP^{k,p}(\g,\mu)=W^{k,p}(\g,\mu)$ (see \cite[Theorem 6.1]{APRR}),
with the advantage that we know explicitly how are
the functions in $W^{k,p}(\g,\mu)$.
Furthermore, we have powerful tools in $W^{k,p}(\g,\mu)$,
like Theorems \ref{t:4.1} and \ref{t:8.3}).

\smallskip

Since, for the sake of generality, we allow
$\|\cdot\|_{W^{k,p}(\g,\mu)}$ to be a seminorm,
it is natural to introduce the following concept.

\begin{definition}
\label{d:K}
Let us consider $1\le p< \infty$ and
a $p$-admissible vectorial measure $\mu$ in $\g$.
Let us define the space $\K (\g, \mu)$ as
$$
\K (\g,\mu) :=\Big\{g:\O^{(0)}\longrightarrow \CC /\,\,
g\in V^{k,p}\big(\g,\mu|_{\O^{(0)}}\big)
,\, \|g\|_{W^{k,p}(\g,\,\mu|_{\O^{(0)}})}=0 \Big\}\,.
$$
\end{definition}
$\K (\g, \mu)$ is the equivalence class of $0$ in
$W^{k,p}(\g,\,\mu|_{\O^{(0)}})$.
Therefore, $\|\cdot\|_{W^{k,p}(\g,\,\mu|_{\O^{(0)}})}$ is a norm
if and only if $\K (\g,\mu)=0$.
It plays an important role in the study of the multiplication operator in Sobolev spaces
(see \cite{RARP2}, \cite{APRR} and Theorem \ref{mult1} below) and in the following
definition, which is crucial in the study of Sobolev spaces
(see \cite{RARP1}, \cite{RARP2} and \cite{APRR}).

In general, it is easy to compute $\K (\g,\mu)$; in the example after
Definition \ref{heuristic}, we can check that
$\K (\g,\mu)=\span \{x, (x_+)^2\}$.

We consider now special classes that we call
$\C_0$ and $\C$. The conditions $(\g,\mu)\in{\C}_0$
and $(\g,\mu)\in\C$ are not very restrictive. The first
one consists, roughly speaking, in considering measures $\mu$ such
that $\| \cdot \|_{W^{k,p}(M_n, \mu)}$ is a norm for some sequence
of compact sets $\{M_n\}$ growing to $\g$. As to the class
$\C$, it is a slight modification of ${\C}_0$, in which we
consider measures $\mu =(\mu_0, \dots , \mu_k)$ such that by
adding a ``minimal" amount of Dirac's deltas to $\mu_0$ we obtain  a measure
in the class ${\C}_0$.

\begin{definition}
\label{d:claseC}
Let us consider $1\le p < \infty$ and
a $p$-admissible vectorial measure $\mu$ in $\g$.
We say that $(\g,\mu)$ belongs
to the class $\C_0$ if there exist compact sets
$M_n$, which are a finite union of compact arcs in $\g$, such that

i$)$ $M_n$ intersects at most a finite number of connected components of
$\O_1\cup\dots\cup\O_k$,

ii$)$ $\K (M_n,\mu)=0$,

iii$)$ $M_n\subseteq M_{n+1}$,

iv$)$ $\cup_n M_{n}=\O^{(0)}$.

\noindent
We say that $(\g,\mu)$ belongs
to the class $\C$ if there exists a measure
$\mu_0'=\mu_0+\sum_{m \in D} c_m\d_{z_m}$ with  $c_m >0$,
$\{z_m\}\subset \O^{(0)}$, $D \subseteq \NN$ and
$(\g,\mu')\in\C_0$,
where $\mu'=(\mu_0',\mu_1,\dots,\mu_k)$ is minimal in the following sense:
there exists $\{M_n\}$ corresponding to $(\g,\mu')\in\C_0$
such that if
$\mu_0''=\mu_0'-c_{m_0}\d_{z_{m_0}}$ with $m_0 \in D$ and
$\mu''=(\mu_0'',\mu_1,\dots,\mu_k)$, then $\K (M_n,\mu'')\neq0$ if
$z_{m_0}\in M_n$.
\end{definition}

\noindent
{\bf Remarks.} \rm

{\bf 1.} The proof in \cite{APRR} of Theorem \ref{t:4.1} below gives that
if $\O^{(0)}\setminus(\O_1\cup\dots\cup\O_k)$
has just a finite number of points in each connected component of $\O^{(0)}$,
and $\K(\g,\mu)=0$, then $(\g,\mu)\in\C_0$.

{\bf 2.} Since the restriction of a function of $\K (\g,\mu)$ to
$M_n$  is in $\K (M_n,\mu)$ for every $n$, then
$(\g,\mu) \in \C_0$ implies $\K(\g,\mu)=0$.

{\bf 3.} If $ (\g,\mu)\in \C_0$, then
$(\g,\mu)\in \C$, with $\mu' =\mu$.

{\bf 4.} The proof in \cite{APRR} of Theorem \ref{t:4.1} below
gives that if for every connected component
$\L$ of $\O_1\cup\cdots\cup\O_k$ we have $\K(\overline\L, \mu)=0$,
then $(\g,\mu)\in\C_0$.
Condition $\#\supp \mu_0|_{\overline\L\cap\O^{(0)}}\ge k$
implies $\K(\overline\L, \mu)=0$.

\begin{teo}
\label{t:5.1}
$($\cite[Theorem 5.1]{APRR}$)$
Let us  consider $1\le p < \infty$ and a $p$-admissible vectorial measure
$\mu=(\mu_0,\dots,\mu_k)$ in a curve $\g$ with $(\g,\mu) \in \C$.
Then the Sobolev space $W^{k,p}(\g,\mu)$ is complete.
\end{teo}

The main ingredient of the proof of this result is Theorem
\ref{t:4.1} below. It allows us to control the $L^\infty$ norm (in
appropriate sets) of a function and its derivatives in terms of
its Sobolev norm. It is also useful by its applications in the
papers \cite{RARP2}, \cite{R1}, \cite{R2}, \cite{R3}, \cite{APRR},
\cite{RY} and \cite{PQRT3}. Furthermore, it is important by itself,
since it answers the following main question: when the
point evaluation functional of $f$ (or $f^{(j)}$) is a
bounded operator in $W^{k,p}(\g,\mu)$?

\begin{teo}
\label{t:4.1}
$($\cite[Theorem 4.1]{APRR}$)$
Let us consider $1\le p < \infty$ and a $p$-admissible
vectorial measure $\mu=(\mu_0,\dots,\mu_k)$ in a curve $\g$.
Let $K_j$ be a finite union of compact arcs contained in  $\O^{(j)}$, for
$0\le j< k$ and $\overline{\mu}$ a right (or left) completion of $\mu$.

{\rm (a)}
If $(\g,\mu)\in\C_0$, then
there exist positive constants $c_1=c_1(\mu, K_0,\dots,K_{k-1})$ and
$c_2 = c_2(\overline{\mu},K_0,\dots,K_{k-1})$ such that
$$
c_1 \sum\limits_{j=0}^{k-1} \|g^{(j)} \| _{L^{\infty}(K_j)} \le
\|g\|_{W^{k,p}(\g,\,\mu)},
\qquad c_2 \, \|g\|_{W^{k,p}(\g,\,\overline{\mu})} \le
\|g\|_{W^{k,p}(\g,\,\mu)},
\qquad \forall \,g \,\in V^{k,p}(\g,\mu).
$$
\indent
{\rm (b)} If $\big(\g,\mu\big)\in \C$
there exist positive constants $c_3=c_3(\mu, K_0,\dots,K_{k-1})$ and
$c_4 = c_4(\overline\mu,K_0,\dots,K_{k-1})$
such that for every $ g\in V^{k,p} (\g,\mu)$, there
exists $g_0 \in V^{k,p}(\g,\mu)$, independent of $K_0, \dots,
K_{k-1}$, $c_3$, $c_4$ and $\overline\mu$, with
$$
\|g_0 -g\|_{W^{k,p} (\g,\,\mu)} =0 \, ,
$$
$$
c_3 \sum\limits_{j=0}^{k-1} \|g_0^{(j)} \| _{L^{\infty}(K_j)} \le
\|g_0\|_{W^{k,p}(\g,\,\mu)}=
\|g\|_{W^{k,p}(\g,\,\mu)}, \qquad c_4  \,
\|g_0\|_{W^{k,p}(\g,\,\overline\mu)} \le
\|g\|_{W^{k,p}(\g,\,\mu)}.
$$
Furthermore, if $g_0,f_0$ are, respectively, these representatives of $g,f$,
we have with the same constants $c_3,c_4$
$$
c_3 \sum\limits_{j=0}^{k-1} \|g_0^{(j)}-f_0^{(j)} \| _{L^{\infty}(K_j)} \le
\|g-f\|_{W^{k,p}(\g,\,\mu)},
 \qquad c_4 \, \|g_0-f_0\|_{W^{k,p}(\g,\,\overline\mu)} \le
\|g - f\|_{W^{k,p}(\g,\,\mu)}.
$$
\end{teo}

\

\section{Basic results on weighted Sobolev spaces.}

In many results about weighted Sobolev spaces in
\cite{RARP1}, \cite{RARP2}, \cite{R1}, \cite{R2}, \cite{R3}, \cite{APRR} and \cite{RY}
appear the hypotheses $(\g,\mu) \in \C$ or $(\g,\mu) \in \C_0$
(see e.g. Theorems \ref{t:5.1} and \ref{t:4.1} above).
The hypotheses $(\g,\mu) \in \C$ and $(\g,\mu) \in \C_0$
appear for the first time in \cite{RARP1},
where they are the sharp conditions in order to obtain
the basic properties of the weighted Sobolev spaces.
Since these conditions are very technical,
it is desirable to obtain some simplifications on them.
Theorem \ref{t:c}, which is the key result in this section, shows a surprising result:
every measure in every curve satisfies $(\g,\mu) \in \C$, and
it has interesting consequences, like Theorems
\ref{t:c1} and \ref{t:c2}.
Theorem \ref{t:c1} is a basic fact in the theory of Sobolev spaces;
it says that $W^{k,p}(\mu)$ is a Banach space for every
$p$-admissible measure $\mu$.

However, the situation with the class $\C_0$ is more difficult.
Theorems \ref{t:cc0}, \ref{t:cc02} and \ref{t:cc03} are simplifications
of the condition $(\g,\mu) \in \C_0$
(we have an example which shows that Theorem \ref{t:cc0} is sharp).
In fact, weighted Sobolev spaces are $W^{1,p}(\g,\mu)$,
at least in ninety per cent of the situations;
in the other cases, they are mostly $W^{2,p}(\g,\mu)$.
Theorem \ref{t:cc02} below guarantees that, in both cases,
$(\g,\mu) \in \C_0$ if and only if $\K(\g,\mu) = 0$,
which is a very simpler condition
($\K (\g,\mu)=0$ means that
$\|\cdot\|_{W^{k,p}(\g,\,\mu|_{\O^{(0)}})}$ is a norm).

\smallskip

We begin with a technical result.

\begin{prop}
\label{p:fj}
Let us  consider $1\le p < \infty$ and a $p$-admissible vectorial measure
$\mu=(\mu_0,\dots,\mu_k)$ in a curve $\g$.
If a connected component $A$ of
$\O_1 \cup \cdots \cup \O_k$ intersects $\O_i$ for some $0 \le i \le k$,
then every function $f \in \K(\g,\mu)$ verifies
$f|_A\in \PP_{i-1}$.
\end{prop}

\noindent
{\bf Remark.}
We use the convention $\PP_{-1} = 0$.

\smallskip

The proof of this result is similar to the proof of \cite[Proposition 2.1]{R4},
which is a version of Proposition \ref{p:fj} for intervals instead of curves.

\begin{teo}
\label{t:c}
Let us  consider any $1\le p < \infty$.
Every $p$-admissible vectorial measure
$\mu=(\mu_0,\dots,\mu_k)$ in any curve $\g$ verifies $(\g,\mu) \in \C$.
\end{teo}

Theorems \ref{t:c}, \ref{t:5.1} and \ref{t:4.1}
have the two following interesting and direct consequences.

\begin{teo}
\label{t:c1}
Let us  consider $1\le p < \infty$ and a $p$-admissible
vectorial measure $\mu=(\mu_0,\dots,\mu_k)$ in a curve $\g$.
Then the Sobolev space $W^{k,p}(\g,\mu)$ is complete.
\end{teo}

\begin{teo}
\label{t:c2}
Let us consider $1\le p < \infty$ and a $p$-admissible
vectorial measure $\mu=(\mu_0,\dots,\mu_k)$ in a curve $\g$.
Let $K_j$ be a finite union of compact arcs contained in  $\O^{(j)}$, for
$0\le j< k$ and $\overline{\mu}$ a right (or left) completion of $\mu$.
Then there exist positive constants $c_3=c_3(\mu, K_0,\dots,K_{k-1})$ and
$c_4 = c_4(\overline\mu,K_0,\dots,K_{k-1})$
such that for every $ g\in V^{k,p} (\g,\mu)$, there
exists $g_0 \in V^{k,p}(\g,\mu)$, independent of $K_0, \dots,
K_{k-1}$, $c_3$, $c_4$ and $\overline\mu$, with
$$
\|g_0 -g\|_{W^{k,p} (\g,\,\mu)} =0 \, ,
$$
$$
c_3 \sum\limits_{j=0}^{k-1} \|g_0^{(j)} \| _{L^{\infty}(K_j)} \le
\|g_0\|_{W^{k,p}(\g,\,\mu)}=
\|g\|_{W^{k,p}(\g,\,\mu)}, \qquad c_4  \,
\|g_0\|_{W^{k,p}(\g,\,\overline\mu)} \le
\|g\|_{W^{k,p}(\g,\,\mu)}.
$$
Furthermore, if $g_0,f_0$ are, respectively, these representatives of $g,f$,
we have with the same constants $c_3,c_4$
$$
c_3 \sum\limits_{j=0}^{k-1} \|g_0^{(j)}-f_0^{(j)} \| _{L^{\infty}(K_j)} \le
\|g-f\|_{W^{k,p}(\g,\,\mu)},
 \qquad c_4 \, \|g_0-f_0\|_{W^{k,p}(\g,\,\overline\mu)} \le
\|g - f\|_{W^{k,p}(\g,\,\mu)}.
$$
\end{teo}

\smallskip

\begin{proof}
In order to prove Theorem \ref{t:c},
without loss of generality we can assume that $\O^{(0)}$ is connected:
if it has connected components $\{O^m\}_m$ with compact sets
$\{M_n^m\}_n$ for each $O^m$ which guarantee that $(\g\cap O^m,\mu) \in \C$, the
compact sets obtained by a diagonal process
$$
M_1^1, M_1^2, M_2^1, M_1^3, M_2^2, M_3^1,
M_1^4, M_2^3, M_3^2, M_4^1, M_1^5, \dots
$$
guarantee that $(\g,\mu) \in \C$
(since there is no relation between the values of a function $f \in W^{k,p}(\g,\mu)$
in two different connected components of $\O^{(0)}$).

Let us consider the connected components $\{\L_\l\}_{\l\in A}$ of
$\O_1\cup\dots\cup\,\O_k$. Recall that
$\O^{(0)}\setminus(\O_1\cup\dots\cup\,\O_k)$ is a discrete set
(see the Remark before Definition \ref{8}).
Moreover this set cannot have any accumulation point in $\O^{(0)}$.
Since $\O^{(0)}$ is connected,
we can take the set of indices $A$ as one of the
following sets:
$\ZZ,\,\ZZ^+,\,\ZZ^-$ or $\{1,\dots,N\}$ for some natural number $N$,
with the property that $\sup\L_\l=\inf\L_{\l+1}=:\b_\l$ if $\l,\l+1\in A$.
Then, $\L_\l=(\b_{\l-1},\b_\l)$.

\smallskip

We need to make some remarks about the structure of $\K(\g,\mu)$.

By Proposition \ref{p:fj},
the functions $f \in \K(\g,\mu)$ verify
$f|_{\L_\l} \in \PP_{j^\l-1}$ in each connected component $\L_\l$,
where $j^\l=\min\{m>0:\, \O_m\cap \L_\l \neq \emptyset\}$.
If we define
$$
j_\l := \min \big\{0 \le m\le j^\l:\,
\sharp \big( \supp \mu_m \cap \O^{(m)}\cap \overline{\L_\l} \big) \ge j^\l -m \big\} \,,
$$
then $f|_{\L_\l} \in \PP_{j_\l-1}$
(with the convention $\PP_{-1}=0$):
we just need to remark that if
$\big( \supp \mu_m \cap \O^{(m)}\cap \overline{\L_\l} \big) \ge j^\l -m$,
then $f^{(m)}|_{\L_\l}$ is zero in at least $j^\l -m$ points;
since $f^{(m)}|_{\L_\l} \in \PP_{j_\l-1-m}$,
then $f^{(m)}|_{\L_\l}=0$ and $f|_{\L_\l} \in \PP_{m-1}$.

Consequently, we can write
$$
f|_{\L_\l}(z)
= b_{j_\l-1}^{\l} z^{j_\l-1} + b_{j_\l-2}^{\l} z^{j_\l-1} + \cdots
+ b_{1}^{\l} z + b_0^{\l} \,.
$$
Since $\supp \mu_j \cap \O^{(j)}\cap \overline{\L_\l}$ is a finite set for
$0\le j < j_\l$,
the other conditions on $f$ in $\L_\l$ can be written as
\begin{equation}
\label{ecuacion11}
\qquad f^{(j)}(z_{ij}^\l ) =0 \qquad
\text{  for  }  0\le j < j_\l, \,\, 1 \le i \le m_j^{\l}\,,
\end{equation}
where the points $\{z_{ij}^\l\}_{ij}$ are
the points with positive $\mu_j$-measure in $\O^{(j)}\cap \overline{\L_\l}$.

Then $(\ref{ecuacion11})$ is a
homogeneous linear system of $m_0^{\l}+ m_{1}^{\l}+ \cdots + m_{j_\l-1}^{\l}$  equations
with the $j_\l$ unknowns $b_{j_\l-1}^{\l}, b_{j_\l-2}^{\l}, \dots, b_{1}^{\l},
b_0^{\l},$ whose solution represents the restriction of the functions in
$\K(\g,\mu)$ to $\L_\l$ in the basis $\{z^{j_\l-1}, z^{j_\l-1}, \dots, z, 1\}$ of
$\PP_{j_\l-1}$.

On the other hand, if $\l,\l+1\in A$, we have
\begin{equation}
\label{ecuacion12}
f^{(j)}(\b_\l^-)=f^{(j)}(\b_\l^+)\,,\qquad \text{ if $\b_\l$ is $j$-regular, }
0\le j < k \,,
\end{equation}
where as usual $f^{(j)}(\b_\l^-)$ and $f^{(j)}(\b_\l^+)$ denotes respectively the left
and right limits of $f^{(j)}$ in $\b_\l$.
Note that we always have $f(\b_\l^-)=f(\b_\l^+)$, since $\b_\l\in \O^{(0)}$
(recall that $\O^{(0)}$ is connected).

Hence, the space $\K(\g,\mu)$ is the solution of the
linear system given by $(\ref{ecuacion11})$ for every $\l \in A$ and
$(\ref{ecuacion12})$ for every $\l$ with $\l,\l+1 \in A$.
In a similar way, the space $\K([\b_{\l^1},\b_{\l^2}],\mu)$ is the solution of
the linear system given by $(\ref{ecuacion11})$ for every $\l^1 < \l \le \l^2$ and
$(\ref{ecuacion12})$ for every $\l^1 < \l < \l^2$.

The coefficient matrix of the system for $\K(\g,\mu)$ have no clear structure:
the equations $(\ref{ecuacion11})$,
for any fixed $\l \in A$,
separate the $j_\l$ unknowns $b_{j_\l-1}^{\l}, b_{j_\l-2}^{\l}, \dots, b_{1}^{\l},
b_0^{\l},$ from the rest of unknowns;
however, the equations $(\ref{ecuacion12})$,
mixed these unknowns with others.
Furthermore, we can not assure that the coefficient matrix represent a bounded
operator in $l^2$ (or even in $l^q$, for some $1 \le q < \infty$);
besides, $\K(\g,\mu)$ is the set of all solutions of these equations
(it is not the set of solutions in some appropriate space).

\medskip

We begin now with the proof.
We assume that the set of indices $A$ is $\ZZ^+$:
if $A=\ZZ^-$ the argument is similar;
if $A=\ZZ$ we just need to combine the argument to the right and to the left of $0$;
if $A$ is a finite set, the proof is direct since then
the equations which define $\K(\g,\mu)$
are a finite linear system.

We assume first that $\b_0 \in \O^{(0)}$.
Let us consider the sequence $\{\dim \K([\b_{0},\b_{\l}],\mu)\}_\l$,
and define
$m := \liminf_{\l \to \infty} \dim \K([\b_{0},\b_{\l}],\mu)$.

\smallskip

If $m < \infty$, there exists a subsequence $\{\l_n\}_n$
with $\dim \K([\b_{0},\b_{\l_n}],\mu) = m$ for every $n$.
We prove first that there exists $m$ points
$z_{1},\dots,z_{m}\in [\b_{0},\b_{\l_1}],$
such that the linear system with the equations of
$\K([\b_{0},\b_{\l_1}],\mu)$
and the $m$ equations $f(z_{1}) = 0,\dots,f(z_{m}) = 0$
has just the trivial solution:

Let us choose a countable dense set $\{q_n\}_n$ in $[\b_{0},\b_{\l_1}]$.
Since any $f\in W^{k,p}(\g,\mu)$ is continuous in $\O^{(0)}$,
the countable system of equations $f(q_n) = 0$
has just the trivial solution $f = 0$ in $\K([\b_{0},\b_{\l_1}],\mu)$.
Hence, if the linear system for $\K([\b_{0},\b_{\l_1}],\mu)$
has $r$ unknowns, we can choose $r$ points
$\{q_{n_1}, \dots, q_{n_r}\}$ such that
the system of linearly independent equations
$f(q_{n_1}) = 0, \dots, f(q_{n_r}) = 0$
has just the trivial solution $f = 0$ in $\K([\b_{0},\b_{\l_1}],\mu)$.
Since $\dim \K([\b_{0},\b_{\l_1}],\mu) = m$,
there exist $z_{1},\dots,z_{m}\in \{q_{n_1}, \dots, q_{n_r}\}$,
such that the linear system given by the equations of $\K([\b_{0},\b_{\l_1}],\mu)$
and $f(z_{1}) = 0,\dots,f(z_{m}) = 0$
has just the trivial solution.
Hence, each equation $f(z_{i}) = 0$ ($1 \le i \le m$)
is linearly independent of the equations of $\K([\b_{0},\b_{\l_1}],\mu)$.
If we define $\mu'=(\mu_0',\mu_1,\dots,\mu_k)$,
where
$\mu_0'=\mu_0+\sum_{i=1}^m \d_{z_i}$,
then $\K([\b_{0},\b_{\l_1}],\mu') = 0$,
and $\mu_0'$ is minimal in the following sense:
if $\mu_0''=\mu_0'-\d_{z_{i_0}}$ with $1 \le i_0 \le m$ and
$\mu''=(\mu_0'',\mu_1,\dots,\mu_k)$,
then $\K ([\b_{0},\b_{\l_1}],\mu'')\neq0$
(in fact, $\dim \K([\b_{0},\b_{\l_1}],\mu'') = 1$).

We prove now that the $m$ equations
$f(z_{1}) = 0,\dots,f(z_{m}) = 0$
are linearly independent with the equations in
$\K ([\b_{0},\b_{\l_n}],\mu)$ for every $n\ge 1$.
Since $\dim \K([\b_{0},\b_{\l_n}],\mu) = m$ for every $n\ge 1$,
this fact proves that $\K([\b_{0},\b_{\l_n}],\mu') = 0$
and $\mu_0'$ is minimal for every $n\ge 1$.
This gives that $(\g,\mu')\in \C_0$ and $(\g,\mu)\in \C$.
We have proved the case $n=1$; therefore, we can consider the case $n>1$.

The unknowns in $\K([\b_{0},\b_{\l_1}],\mu)$
appear also in the $m$ equations
$f(z_{1}) = 0,\dots,f(z_{m}) = 0$,
and in
\begin{equation}
\label{ecuacion13}
f^{(j)}(\b_{\l_1}^-)=f^{(j)}(\b_{\l_1}^+)\,,\qquad \text{ if $\b_{\l_1}$ is $j$-regular, }
0\le j < k \,.
\end{equation}
The equations of $\K([\b_{0},\b_{\l_1}],\mu)$ and
$f(z_{1}) = 0,\dots,f(z_{m}) = 0$
are the equations of $\K([\b_{0},\b_{\l_1}],\mu')$.
We know that $\K([\b_{0},\b_{\l_1}],\mu') = 0$,
and then these equations have just the trivial solution.
Hence, $f=0$ in $[\b_{0},\b_{\l_1}]\cap \O^{(0)}$, and
$(\ref{ecuacion13})$ is equivalent to
\begin{equation}
\label{ecuacion14}
0=f^{(j)}(\b_{\l_1}^+)\,,\qquad \text{ if $\b_{\l_1}$ is $j$-regular, }
0\le j < k \,.
\end{equation}
Since the equations in $\K([\b_{0},\b_{\l_n}],\mu')$
are equivalent to the equations in $\K([\b_{0},\b_{\l_1}],\mu')$,
$(\ref{ecuacion14})$ and the equations in $\K([\b_{\l_1},\b_{\l_n}],\mu)$,
and the unknowns in $\K([\b_{0},\b_{\l_1}],\mu')$
are neither in
$(\ref{ecuacion14})$ nor $\K([\b_{\l_1},\b_{\l_n}],\mu)$,
the equations in $\K([\b_{0},\b_{\l_1}],\mu')$
(in particular, $f(z_{1}) = 0,\dots,f(z_{m}) = 0$)
are linearly independent with the equations
$(\ref{ecuacion14})$ and $\K([\b_{\l_1},\b_{\l_n}],\mu)$.
Since we have proved above that
$f(z_{1}) = 0,\dots,f(z_{m}) = 0$
are linearly independent with the equations of
$\K ([\b_{0},\b_{\l_1}],\mu)$, we deduce that
$f(z_{1}) = 0,\dots,f(z_{m}) = 0$
are linearly independent with the equations of
$\K ([\b_{0},\b_{\l_n}],\mu)$.
Since $\dim \K([\b_{0},\b_{\l_n}],\mu) = m$ for every $n\ge 1$,
this fact proves that $\K([\b_{0},\b_{\l_n}],\mu') = 0$
and $\mu_0'$ is minimal for every $n\ge 1$.
Then, $(\g,\mu')\in \C_0$ and $(\g,\mu)\in \C$.

\smallskip

If $m = \infty$, consider $\l_1$ with
$\dim \K([\b_{0},\b_{\l_1}],\mu)= \min \{\dim \K([\b_{0},\b_{\l}],\mu): \, \l\in \ZZ^+
\}:= m_1$. As in the last case, we can find
$z_{1},\dots,z_{m_1}\subset [\b_{0},\b_{\l_1}]$,
such that the linear system given by the equations of $\K([\b_{0},\b_{\l_1}],\mu)$
and $f(z_{1}) = 0,\dots,f(z_{m_1}) = 0$
has just the trivial solution.
Hence, each equation $f(z_{i}) = 0$ ($1 \le i \le m_1$)
is linearly independent of the equations of $\K([\b_{0},\b_{\l_1}],\mu)$.
If we define $\mu^1=(\mu_0^1,\mu_1,\dots,\mu_k)$, where
$\mu_0^1=\mu_0+\sum_{i=1}^{m_1} \d_{z_i}$,
then $\K([\b_{0},\b_{\l_1}],\mu^1) = 0$,
and $\mu_0^1$ is minimal.

Let us assume now that we have $\l_r$ and
$\mu^r=(\mu_0^r,\mu_1,\dots,\mu_k)$, where
$\mu_0^r=\mu_0+\sum_{i=1}^{m_r} \d_{z_i}$,
$\K([\b_{0},\b_{\l_r}],\mu^r) = 0$,
and $\mu_0^r$ is minimal,
for $1 \le r \le n$.
We choose $\l_{n+1}$ with
$\dim \K([\b_{0},\b_{\l_{n+1}}],\mu)= \min \{\dim
\K([\b_{0},\b_{\l}],\mu): \, \l> \l_n \}:= m_{n+1}-m_n$.
As above, we can find
$z_{m_n+1},\dots,z_{m_{n+1}}\subset [\b_{0},\b_{\l_{n+1}}]$,
such that the linear system given by the equations of $\K([\b_{0},\b_{\l_{n+1}}],\mu^n)$
and $f(z_{m_n+1}) = 0,\dots,f(z_{m_{n+1}}) = 0$
has just the trivial solution.
Since $\K([\b_{0},\b_{\l_{n}}],\mu^n) = 0$,
it is clear that
$z_{m_n+1},\dots,z_{m_{n+1}}\subset (\b_{\l_n},\b_{\l_{n+1}}]$,
and each equation $f(z_{i}) = 0$ ($m_{n}+1 \le i \le m_{n+1}$)
is linearly independent of the equations of $\K([\b_{0},\b_{\l_{n+1}}],\mu^n)$.
If we define $\mu^{n+1}=(\mu_0^{n+1},\mu_1,\dots,\mu_k)$, where
$\mu_0^{n+1}=\mu_0+\sum_{i=1}^{m_{n+1}} \d_{z_i}$,
then $\K([\b_{0},\b_{\l_{n+1}}],\mu^{n+1}) = 0$,
and $\mu_0^{n+1}$ is minimal.
This process defines inductively the measures $\{\mu^{n}\}_n$ and the points
$\{z_i\}_{i=1}^{l}$,
with $l:= \lim_{n\to \infty} m_n$,
and we can consider $\mu'$ with
$\mu_0':=\mu_0+\sum_{i=1}^{l} \d_{z_i}$,
and $M_n := [\b_{0},\b_{\l_{n}}]$.
The argument as above shows that
$(\g,\mu')\in \C_0$ and $(\g,\mu)\in \C$.

\smallskip

If $\b_0 \notin \O^{(0)}$,
we can not use the same construction since
$[\b_0,\b_{\l_n}] \nsubseteq \O^{(0)}$.
We can avoid this problem by choosing a sequence
$\{\a_{n}\}_n \in (\b_0,\b_{1})$ converging to $\b_0$ and such that
$\supp \mu_j \cap \Omega^{(j)} \cap (\b_0,\b_{1}] \subseteq [\a_n,\b_{1}]$
for every $n$ and every $0 \le j < j_{\l_1}$:
it is enough to consider $M_n := [\a_{n},\b_{\l_{n}}]$.
\end{proof}

We deal now with the class $\C_0$.

\begin{teo}
\label{t:cc0}
Let us  consider $1\le p < \infty$ and a $p$-admissible
vectorial measure $\mu=(\mu_0,\dots,\mu_k)$ in a curve $\g$.
Let us assume that every function $f \in \K(\g,\mu)$ verifies
$f''=0$ a.e. in $\O_1 \cup \cdots \cup \O_k$.
Then, $(\g,\mu) \in \C_0$ if and only if $\K(\g,\mu) = 0$.
\end{teo}

One can think that, in a similar way to Theorem \ref{t:c}, perhaps
the conclusion of Theorem \ref{t:cc0} holds without the hypothesis on $f''$.
However, this is not true. In fact,
Theorem \ref{t:cc0} is sharp in the following sense:
its conclusion does not hold if we substitute the hypothesis
``$f''=0$ a.e. in $\O_1 \cup \cdots \cup \O_k$" by
``$f'''=0$ a.e. in $\O_1 \cup \cdots \cup \O_k$", as the following example shows.

\medskip

\noindent
{\bf Example.}
Fix $1\le p < \infty$ and consider
a finite $p$-admissible vectorial measure
$\mu=(\mu_0,\mu_1,\mu_2,\mu_3)$
in the compact interval $\g:=[0,1]$ defined as follows:
$$
\begin{aligned}
\mu_0 & := \sum_{m=0}^{\infty} 2^{-m} \d_{2^{-2m-1}} \,, \qquad
\mu_1:= \sum_{m=0}^{\infty} 2^{-m} \d_{3\cdot 2^{-2m-2}} \,,
\\
\mu_2 & := 0 \,, \qquad
w_3(x):= \sum_{m=0}^{\infty} (2^{-2m} - x)^{2p-1} (x-2^{-2m-2})^{2p-1}
\chi_{_{[2^{-2m-2}, 2^{-2m}]}}(x) \,.
\end{aligned}
$$
It is not difficult to check that
$\O_1=\O_2= \emptyset$,
$\O_3=(0,1) \setminus \cup_{m=1}^{\infty} \{2^{-2m}\}$,
$\O^{(0)}=(0,1]$ (since $2p-1<3p-1$; see Remark $1$(b) after
Definition \ref{d:regular})
and $\O^{(1)}=\O^{(2)}=(0,1) \setminus \cup_{m=1}^{\infty} \{2^{-2m}\}$.

Let us see first that $\K(\g,\mu) = 0$.
Consider $f \in \K(\g,\mu)$. Then $f'''=0$ a.e. in $(2^{-2m-2}, 2^{-2m})$
and hence $f_m:= f|_{(2^{-2m-2}, 2^{-2m})} \in \PP_2$ for each $m \ge 0$.
It is not difficult to see that the equations
$f_m(2^{-2m-1})=0$, $f_m'(3\cdot 2^{-2m-2})=0$
are equivalent to the equations
$f_m(2^{-2m-1})=0$, $f_m(2^{-2m})=0$
($3\cdot 2^{-2m-2}$ is the middle point between $2^{-2m-1}$ and $2^{-2m}$).
Since $2^{-2m-2} \in \O^{(0)}$, then
$f_{m+1}(2^{-2m-2})=f_m(2^{-2m-2})=0$,
and consequently
$f_m(2^{-2m-2})=0$, $f_m(2^{-2m-1})=0$, $f_m(2^{-2m})=0$, for each $m \ge 0$.
Then, $f_m=0$ for each $m \ge 0$, and $f=0$ in $(0,1]$.
Hence, $\K(\g,\mu) = 0$.

However, $(\g,\mu) \notin \C_0$:
if we define $M_n:=[2^{-2n-2}, 1]$, we can find a function
$g_n \in \K(M_n,\mu) \setminus 0$, for instance,
$$
g_n:=
\begin{cases}
(x - 2^{-2n} ) (x-2^{-2n-1}) \,,
&\qquad \text{if } \, x \in [2^{-2n-2}, 2^{-2n}] \,,
\\
0\,,  &\qquad \text{if } \, x \in [2^{-2n}, 1]\, .
\end{cases}
$$
It is easy to check that for any choice of compact sets
$M_n \subset (0,1]$ we also have
$\K(M_n,\mu) \neq 0$, and then $(\g,\mu) \notin \C_0$.

\medskip

\begin{proof}
As in the proof of Theorem \ref{t:c},
without loss of generality we can assume that $\O^{(0)}$ is connected,
and we can consider the connected components $\{\L_\l\}_{\l\in A}$ of
$\O_1\cup\dots\cup\,\O_k$. Recall that
we can take the set of indices $A$ as one of the
following sets:
$\ZZ,\,\ZZ^+,\,\ZZ^-$ or $\{1,\dots,N\}$ for some natural number $N$,
with the property that $\sup\L_\l=\inf\L_{\l+1}$ if $\l,\l+1\in A$.

\smallskip

We need to make some remarks about the structure of $\K(\g,\mu)$
in this particular case.

The functions $f \in \K(\g,\mu)$ verify
$f''=0$ a.e. in every connected component $\L_\l$;
then $f'|_{\L_\l}$ is constant, since $f' \in AC_{loc}(\O^{(1)})
\subseteq AC_{loc}(\L_\l)$, and consequently $f|_{\L_\l} \in \PP_1$;
furthermore, the other conditions on $f$ in $\L_\l$ can be written as
\begin{equation}
\label{ecuacion1}
\qquad f^{(j)}(z_j^i ) =0 \qquad
\text{  for  }  j =0,1, \,\, 1 \le i \le m_j\,.
\end{equation}
\indent
Roughly speaking, the points $\{z_j^i\}_{ij}$ are
the points with positive $\mu_j$-measure in $\L_\l$.
More concretely:
If $\mu_1\big(\overline{\L_\l} \cap \O^{(1)}\big)>0$,
then $f'|_{\L_\l}=0$ and $f|_{\L_\l} \in \PP_0$;
this condition is equivalent to
$f'(z_1^1)=0$ for any $z_1^1 \in \L_\l$.
If $\supp \mu_0$ has at least two points $\{z_0^1,z_0^2\}$
in $\overline{\L_\l} \cap \O^{(0)}$,
then $f|_{\L_\l}=0$ and this condition is equivalent to
$f(z_0^1)=f(z_0^2)=0$;
if $\supp \mu_0$ has exactly one point $\{z_0^1\}$
in $\overline{\L_\l} \cap \O^{(0)}$,
then the condition is $f(z_0^1)=0$;
if $\mu_0\big(\overline{\L_\l} \cap \O^{(0)}\big)=0$,
then there is no additional condition.

If we write in $\L_\l$
$$
f(z)= \d_{1} z + \d_0\,,
$$
we see that $(\ref{ecuacion1})$ is a
homogeneous linear system of $m_0+ m_{1}$  equations
with the two unknowns $\d_{1}, \d_0$,
whose solution represents the restriction of the functions in
$\K(\g,\mu)$ to $\L_\l$ in the basis
$\{z, 1\}$ of $\PP_{1}$.
Note that the constants $\d_i$ and $m_i$ obviously depend on $\l$.

On the other hand, if $\l,\l+1\in A$ and $\b_\l:=\sup \L_\l=\inf \L_{\l+1}$,
we have
\begin{equation}
\label{ecuacion2}
f^{(j)}(\b_\l^-)=f^{(j)}(\b_\l^+)\,,\qquad \text{ if $\b_\l$ is $j$-regular, }  j=0,1\,,
\end{equation}
where as usual $f^{(j)}(\b_\l^-)$ and $f^{(j)}(\b_\l^+)$ denotes respectively the left
and right limits of $f^{(j)}$ in $\b_\l$.
Note that we always have $f(\b_\l^-)=f(\b_\l^+)$, since $\b_\l\in \O^{(0)}$
(recall that $\O^{(0)}$ is connected).

Let us remark that the condition $f'(\b_\l^-)=f'(\b_\l^+)$
(which appears if and only if $\b_\l\in \O^{(1)}$) implies that
$f|_{\L_\l\cup \{\b_\l\} \cup \L_{\l+1}} \in \PP_1$.
Consequently, it is natural to consider the connected components $\{\G_m\}_{m\in B}$ of
$\O_1\cup\dots\cup\,\O_k \cup \{\b_\l:\,\b_\l \in \O^{(1)}\}$;
we can take the set of indices $B$ as one of the following sets:
$\ZZ,\,\ZZ^+,\,\ZZ^-$ or $\{1,\dots,N\}$ for some natural number $N$,
with the property that $\G_m=(a_{m-1},a_m)$
(then $\sup\G_m=\inf\G_{m+1}=a_m$ if $m,m+1\in B$).
We have $f|_{\G_m} \in \PP_1$ for each $m \in B$.
The other conditions on $f$ in $\G_m$ can be written as
\begin{equation}
\label{ecuacion3}
\qquad f^{(j)}(\z_j^i ) =0 \qquad
\text{  for  }  j =0,1, \,\, 1 \le i \le n_j\,.
\end{equation}
If we write in $\G_m$
$$
f(z)= \alpha_{1} z + \alpha_0\,,
$$
we see that $(\ref{ecuacion3})$ is a
homogeneous linear system of $n_0+ n_{1}$  equations
with the two unknowns $\alpha_{1}, \alpha_0$,
whose solution represents the restriction of the functions in
$\K(\g,\mu)$ to $\G_m$ in the basis
$\{z, 1\}$ of $\PP_{1}$.
Note that the constants $\a_i$ and $n_i$ obviously depend on $m$.

On the other hand, we have
\begin{equation}
\label{ecuacion4}
f(a_m^-)=f(a_m^+)\,,\qquad \text{ if }  m,m+1 \in B\,.
\end{equation}
Then we have that
$\K(\g,\mu)$ is the solution of the
linear system given by $(\ref{ecuacion3})$ for every $m\in B$ and $(\ref{ecuacion4})$
for every $m\in B$ such that $m+1\in B$
(we also have that $\K(\g,\mu)$ is the solution of the
linear system given by $(\ref{ecuacion1})$ for every $\l\in A$ and $(\ref{ecuacion2})$
for every $\l\in A$ such that $\l+1\in A$).
Consequently, the elements of $\K(\g,\mu)$ are linear splines.

We describe now an algorithm, which we call \emph{extension process},
in order to construct a function in
$\K \big( \cup_{m\ge m_0} \overline{\G_m} ,\mu \big) \setminus 0$
under the hypothesis $\K \big( \G_m ,\mu \big) \neq 0$ for every $m\ge m_0$.
Given a function $f \in \K \big( \G_m ,\mu \big) \setminus 0$,
we can extend it to a function
$f \in \K \big( \G_m \cup \{a_m\} \cup \G_{m+1} , \mu \big) \setminus 0$
as follows.
If $f(a_m)=0$, it is enough to define $f=0$ in $\G_{m+1}$.
If $f(a_m)=\d \neq 0$, we need to define $f$ in $\G_{m+1}$ as
$$
f(z)= \alpha_{1} z + \alpha_0 \,,
$$
with $f(a_m)=\d$
and perhaps one additional condition
(it is not possible to have more than one condition, since then
$\K \big( \G_{m+1} ,\mu \big) = 0$,
which is a contradiction).
If there is no more condition, we define $f$ in $\G_{m+1}$ as
the function with graph the straight line joining $(a_m,\d)$
with $(a_{m+1},0)$;
if we have the condition $f(\z_0^1) =0$, we define $f$ in $\G_{m+1}$ as
the function with graph the straight line joining $(a_m,\d)$
with $(\z_0^1,0)$;
if we have the condition $f'(\z_1^1) =0$, we define $f=\d$ in $\G_{m+1}$.
Now, since $\K \big( \G_{m_0} ,\mu \big) \neq 0$ we can choose a function
$f \in \K \big( \G_{m_0} ,\mu \big) \neq 0$;
applying inductively the extension process, we can extend $f$ to a function in
$\K \big( \cup_{m\ge m_0} \overline{\G_m} ,\mu \big) \setminus 0$.

\medskip

We begin now with the proof.
Remark $2$ to Definition \ref{d:claseC} shows that
if $(\g,\mu) \in \C_0$, then $\K(\g,\mu) = 0$.
Assume now that $\K(\g,\mu) = 0$.
Consider the connected components $\{\G_m\}_{m\in B}$ of
$\O_1\cup\dots\cup\,\O_k \cup \{\b_\l:\,\b_\l \in \O^{(1)}\}$.

It is clear that there exists $m\in B$ with
$\K \big( \G_m ,\mu \big) = 0$,
since otherwise the extension process
for $m\ge m_0$ and $m\le m_0$, gives a function
$f \in \K (\g ,\mu ) \setminus 0$,
which is a contradiction.

Consider the ordered subset $\{m_j\}_j \subseteq B$
with $\K \big( \G_{m_j} , \mu \big) = 0$.

$(A)$
Let us note first that if $i<j$, then
$\K \big( \cup_{m=m_i}^{m_j} \overline{\G_m} , \mu \big) = 0$.
It is trivial if $m_j=m_i+1$; then we can assume that $m_j> m_i+1$.
This space is the solution of the
linear system given by $(\ref{ecuacion3})$ for every $m_i \le m \le m_j$ and
$(\ref{ecuacion4})$ for every $m_i \le m < m_j$;
since the unknowns in $\K \big( \G_{m_i} , \mu \big)$ and
$\K \big( \G_{m_j} , \mu \big)$ are $0$, the linear system with
the unknowns in
$\K \big( \cup_{m=m_i+1}^{m_j-1} \overline{\G_m} , \mu \big)$
 is ``isolated"
of the rest of equations of $\K(\g,\mu)$:
conditions
$$
f(a_{m_i}^-)=f(a_{m_i}^+)\,, \qquad
f(a_{m_j-1}^-)=f(a_{m_j-1}^+)\,,
$$
are in fact
$$
0=f(a_{m_i}^+)\,, \qquad
f(a_{m_j-1}^-)=0\,,
$$
since $\K(\G_{m_i},\mu) = 0$ and $\K(\G_{m_j},\mu) = 0$.
Since $\K(\g,\mu) = 0$,
every unknown in
$\K \big( \cup_{m=m_i+1}^{m_j-1} \overline{\G_m} , \mu \big)$
is $0$, and we conclude that
$\K \big( \cup_{m=m_i}^{m_j} \overline{\G_m} , \mu \big) = 0$.

$(B)$
If $\{m_j\}_j$ has a maximum $m^*$, then for any $M > m^*$ we have
$\K \big( \cup_{m= m^*}^M \overline{\G_m} , \mu \big) = 0$:
Seeking for a contradiction assume that there exists
$f \in \K \big( \cup_{m= m^*}^M \overline{\G_m} , \mu \big) \setminus 0$.
The extension process
(since $\K \big( \G_{m} , \mu \big) \neq 0$ for every $m > m^*$)
gives a function
$f \in \K \big( \cup_{m\ge m^*} \overline{\G_m} , \mu \big) \setminus 0$.
Since $\K ( \G_{m^*} , \mu ) = 0$,
then $f=0$ in $\overline{\G_{m^*}}$
and we can extend it as $f=0$ to
$\cup_{m< m^*} \overline{\G_m}$,
obtaining a function in $\K (\g ,\mu ) \setminus 0$,
which is a contradiction.
Then, we conclude that
$\K \big( \cup_{m= m^*}^M \overline{\G_m} , \mu \big) = 0$
for any $M> m^*$.

$(C)$
A similar argument proves that if
$\{m_j\}_j$ has a minimum $m_*$, then for any $M < m_*$ we have
$\K \big( \cup_{m=M}^{m_*} \overline{\G_m} , \mu \big) = 0$.

The facts $(A)$, $(B)$ and $(C)$ allows to choose easily
the sets $M_n$ with $\K (M_n,\mu)=0$ verifying the properties
which guarantee that $(\g,\mu) \in \C_0$:

Assume that $B= \ZZ$ (the other cases are easier).

If $\lim_{j \to -\infty} m_j = - \infty$
and $\lim_{j \to \infty} m_j = \infty$,
we can choose $M_n = \cup_{m= m_{-n}}^{m_n} \overline{\G_m}\,$.

If $\lim_{j \to -\infty} m_j = - \infty$
and $m_j \le m^*$ for every $j\in \ZZ$,
we can choose $M_n = \cup_{m= m_{-n}}^{m^* + n} \overline{\G_m}\,$.

If $m_* \le m_j$ for every $j\in \ZZ$
and $\lim_{j \to \infty} m_j = \infty$,
we can choose $M_n = \cup_{m= m_*-n}^{m_n} \overline{\G_m}\,$.

Finally, if $m_* \le m_j \le m^*$ for every $j\in \ZZ$,
we can choose $M_n = \cup_{m= m_*-n}^{m^*+n} \overline{\G_m}\,$.
\end{proof}

Theorem \ref{t:cc0} and Proposition \ref{p:fj} give the following:

\begin{teo}
\label{t:cc02}
Let us  consider $1\le p < \infty$ and a $p$-admissible
vectorial measure $\mu=(\mu_0,\dots,\mu_k)$ in a curve $\g$.
Let us assume that every connected component of
$\O_1 \cup \cdots \cup \O_k$ intersects $\O_0 \cup \O_1 \cup \O_2$.
Then, $(\g,\mu) \in \C_0$ if and only if $\K(\g,\mu) = 0$.
\end{teo}

As a corollary we obtain the following result.

\begin{teo}
\label{t:cc03}
Let us  consider $1\le p < \infty$ and
a $p$-admissible vectorial measure $\mu=(\mu_0,\dots,\mu_k)$
in a curve $\g$, with $k=1$ or $k=2$.
Then, $(\g,\mu) \in \C_0$ if and only if $\K(\g,\mu) = 0$.
\end{teo}

Theorems \ref{t:cc03} and \ref{t:4.1} give the following direct consequence.

\begin{teo}
Let us consider $1\le p < \infty$ and a $p$-admissible vectorial measure
$\mu=(\mu_0,\dots,\mu_k)$ in a curve $\g$, with $k=1$ or $k=2$.
Let $K_j$ be a finite union of compact arcs contained in  $\O^{(j)}$, for
$0\le j< k$ and $\overline{\mu}$ a right (or left) completion of $\mu$.
If $\K(\g,\mu)= 0$, then
there exist positive constants $c_1=c_1(\mu, K_0,K_{k-1})$ and
$c_2 = c_2(\overline{\mu},K_0,K_{k-1})$ such that
$$
c_1 \sum\limits_{j=0}^{k-1} \|g^{(j)} \| _{L^{\infty}(K_j)} \le
\|g\|_{W^{k,p}(\g,\,\mu)},
\qquad c_2 \, \|g\|_{W^{k,p}(\g,\,\overline{\mu})} \le
\|g\|_{W^{k,p}(\g,\,\mu)},
\qquad \forall \,g \,\in V^{k,p}(\g,\mu).
$$
\end{teo}

\

\section{Results on the multiplication operator.}

Recall that when every polynomial has finite
$W^{k,p}(E,\mu)$-norm, we denote by $\PP^{k,p}(E,\mu)$ the
completion of $\PP$ with that norm.
Since our aim is to bound the multiplication operator in $\PP^{k,p}(E,\mu)$,
in this section we just consider measures such that every polynomial has finite
Sobolev norm. Hence, for any $0 \le j \le k$,
$$
\mu_j(\CC)^{1/p} = \big\|1\big\|_{L^{p}(E,\mu_j)}
\le \big\|z^j/j!\big\|_{W^{k,p}(E,\mu)} < \infty \,,
$$
and consequently, $\mu$ is finite.

M. Castro and A. Dur\'an  \cite{CD}
proved that if the multiplication operator is bounded in $\PP^{k,p}(\mu)$
then the support of $\mu$ is compact.
Then, we just need to consider finite vectorial measures with compact support.
In this case, we usually have $\PP^{k,p}(\g,\mu)=W^{k,p}(\g,\mu)$
(see \cite[Theorem 6.1]{APRR}).

\smallskip

First of all, some remarks about the definition of the
multiplication operator.
We start with a definition which makes sense for measures defined
in arbitrary Borel sets $E\subseteq \CC$ (not necessarily curves).

\begin{definition}
\label{d:8.1}
If $\mu$ is a vectorial measure in the
Borel set $E\subseteq \CC$,
we say that the multiplication operator is \emph{well defined in}
$\PP^{k,p}(E,\mu)$ if given any sequence $\{s_n\}$ of polynomials
converging to $0$ in the $W^{k,p}(E,\mu)$-norm,
then $\{zs_n\}$ also converges to $0$ in the $W^{k,p}(E,\mu)$-norm.
In this case, if $\{q_n\}\in \PP^{k,p}(E,\mu)$, we define
$\M(\{q_n\}):=\{zq_n\}$, where $z$ is the independent variable.
If we choose another Cauchy sequence
$\{r_n\}$ representing the same element in $\PP^{k,p}(E,\mu)$
(i.e. $\{q_n-r_n\}$ converges to $0$ in the $W^{k,p}(E,\mu)$-norm),
then $\{zq_n\}$ and $\{zr_n\}$ represent the same element in $\PP^{k,p}(E,\mu)$
(since $\{z(q_n-r_n)\}$ converges to $0$ in the $W^{k,p}(E,\mu)$-norm).
\end{definition}

We can also think of another definition which is as natural
as the previous one in the case of curves.

\begin{definition}
\label{d:8.2}
If $\mu$ is a $p$-admissible vectorial measure in
$\g$ $($and hence $W^{k,p}(\g,\mu)$
is a space of classes of functions$)$, we say that
the multiplication operator is \emph{well defined in}
$W^{k,p}(\g,\mu)$ if given any function $h\in V^{k,p}(\g,\mu)$
with $\|h\|_{W^{k,p}(\g,\mu)}=0$, we have
$\|zh\|_{W^{k,p}(\g,\mu)}=0$.
In this case, if $[f]$ is an equivalence class in $W^{k,p}(\g,\mu)$,
we define $\M([f]):=[zf]$. If we choose another
representative $g$ of $[f]$ $($i.e. $\|f-g\|_{W^{k,p}(\g,\mu)}=0)$
we have $[zf]=[zg]$, since $\|z(f-g)\|_{W^{k,p}(\g,\mu)}=0$.
\end{definition}

Although both definitions are natural, it is possible for
a $p$-admissible measure $\mu$ with $W^{k,p}(\g,\mu)=\PP^{k,p}(\g,\mu)$
that $\M$ is well defined in $W^{k,p}(\g,\mu)$ and
it is not well defined in $\PP^{k,p}(\g,\mu)$
(see the example after \cite[Theorem 4.2]{R4}).
The following elementary lemma characterizes the spaces $\PP^{k,p}(E,\mu)$
for which $M$ is well defined.

\begin{lemma}
\label{l:8.1}
$($\cite[Lemma 8.1]{APRR}$)$
Let us consider $1\le p<\infty$ and $\mu=(\mu_0,\dots,\mu_k)$ a
vectorial measure in a Borel set $E\subseteq\CC$. The following facts are equivalent.

$(1)$ The multiplication operator is well defined in
$\PP^{k,p}(E,\mu)$.

$(2)$ The multiplication operator is bounded in $\PP^{k,p}(E,\mu)$.

$(3)$ There exists a positive constant $c$ such that
$$
\|z q\|_{W^{k,p}(E,\mu)} \le c\, \|q\|_{W^{k,p}(E,\mu)} \,,
\qquad \text{for every } q\in \PP\,.
$$
\end{lemma}

\begin{definition}
A vectorial  measure $\mu = (\mu_0, \dots , \mu_k)$ in the complex plane
is \emph{extended sequentially dominated} (and we write $\mu \in ESD$)
if there exists a positive constant $c$
such that $\mu_{j+1}\le c \, \mu_{j}$ for $0\le j<k$.
\end{definition}

This kind of measures plays a main role in the study of the multiplication operator:

\begin{teo}
\label{t:8.1}
$($\cite[Theorem 8.1]{APRR}$)$
Let us consider $1\le p<\infty$ and $\mu=(\mu_0,\dots,\mu_k)$ a
finite vectorial measure in a compact set $E$.
Then, the multiplication operator is bounded in $\PP^{k,p}(E,\mu)$
if and only if there exists a vectorial measure $\mu'\in ESD$
such that the Sobolev norms in
$W^{k,p}(E,\mu)$ and $W^{k,p}(E,\mu')$ are comparable on $\PP$.
Furthermore, we can choose $\mu'=(\mu_0',\dots,\mu_k')$
with $\mu_j':=\mu_j+\mu_{j+1}+\cdots+\mu_k$.
\end{teo}

Although this result characterizes the measures with $\M$ bounded,
it is convenient to obtain more practical criteria in order to guarantee
the boundedness of $\M$.

\medskip

If we consider the case of a curve $E=\g$,
we have the following results.

First, let us note that
the multiplication operator $\M$ is bounded in $W^{k,p}(\g,\mu)$
if and only if there exists a positive constant $c$ such that
$$
\|z f\|_{W^{k,p}(\g,\mu_j)} \le c\,
\|f\|_{W^{k,p}(\g,\mu)}        \,,
$$
for every $f\in V^{k,p}(\g,\mu)$.
Consequently, if $\M$ is bounded in $W^{k,p}(\g,\mu)$
and $\PP \subseteq W^{k,p}(\g,\mu)$,
then it is bounded in $\PP^{k,p}(\g,\mu)$,
since $W^{k,p}(\g,\mu)$ is a complete space by Theorem \ref{t:c1}.

\smallskip

The following result characterizes when
$\M$ is a well defined operator in $W^{k,p}(\g,\mu)$.

\begin{teo}
\label{t:8.3}
Let us consider $1\le p<\infty$ and a
$p$-admissible vectorial measure $\mu$ in $\g$.
Then the multiplication operator $\M$ is well defined
in $W^{k,p}(\g,\mu)$ if and only if $\K (\g,\mu)=0$.
\end{teo}

\smallskip

The proof of this result is similar to the proof of \cite[Theorem 4.2]{R4},
which is a version of Theorem \ref{t:8.3} for intervals instead of curves.

\smallskip

One can think that, in a similar way to Lemma \ref{l:8.1},
the multiplication operator $\M$ is well defined
in $W^{k,p}(\g,\mu)$ if and only if it is bounded in $W^{k,p}(\g,\mu)$.
However, this is not true, as the example after \cite[Theorem 4.2]{R4} shows.

\begin{lemma}
\label{l:8.2}
Let us consider $1\le p<\infty$ and a
$p$-admissible vectorial measure $\mu$ in a compact curve $\g$.
Then, the multiplication operator $\M$ is bounded in $W^{k,p}(\g,\mu)$
if and only if there exists a positive constant $c$ such that
$$
\|f^{(j-1)}\|_{L^{p}(\g,\mu_j)} \le c\,
\|f\|_{W^{k,p}(\g,\mu)}        \,,
$$
for every $1\le j\le k$ and $f\in V^{k,p}(\g,\mu)$.
\end{lemma}

\smallskip

The proof of this result is similar to the proof of \cite[Lemma 4.2]{R4},
which is a version of Lemma \ref{l:8.2} for intervals instead of curves.

\begin{definition}
\label{d:consistent}
Let us consider $1 \le p < \infty$, a curve $\g$ and $z_0,z_1 \in \g$.
We say that a weight $w$ in $\g$ is \emph{right-consistent}
(respectively, \emph{left-consistent}) in $[z_0,z_1]$ if
$\L_{p,[z_0,z_1]}^+ (w, w)<\infty$
(respectively, $\L_{p,[z_0,z_1]}^- (w, w)<\infty$).
\end{definition}

\noindent
{\bf Remark.}
It is easy to see that if a weight $w$ is comparable to a
non-decreasing function $w^0$ in $[z_0,z_1]$
(respectively, non-increasing),
then it is right-consistent
(respectively, left-consistent) in $[z_0,z_1]$:
$$
\begin{aligned}
\L_{p,[z_0,z_1]}^+ (w, w)
& = \sup_{z_0 < z < z_1} \Big(\int_{z_0}^{z} w \Big)
\Big( \int_{z}^{z_1} w^{-1/(p-1)} \Big)^{p-1}
\le \, c \!\! \sup_{z_0 < z < z_1} \Big(\int_{z_0}^{z} w^0 \Big)
\Big( \int_{z}^{z_1} (w^0)^{-1/(p-1)} \Big)^{p-1}
\\
& \le \, c \!\! \sup_{z_0 < z < z_1} L([z_0,z]) \, w^0 (z) \,
\Big( L([z,z_1]) \, w^0(z)^{-1/(p-1)} \Big)^{p-1}
\\
& = \, c \!\! \sup_{z_0 < z < z_1} L([z_0,z]) \, L([z,z_1])^{p-1} < \infty \,.
\end{aligned}
$$
Note that $w=0$ in a right-neighborhood of $z_0$ (or in a
left-neighborhood of $z_1$), and even in $[z_0,z_1]$, is allowed.

\begin{lemma}
\label{l:8.3}
Let us consider $1\le p<\infty$, a curve $\g$, $z_0,z_1\in \g$ and a
right-consistent weight $w$ in $[z_0,z_1]$.
Then the largest open set $U\subset (z_0,z_1)$ with $w\in B_p(U)$
is $(\a,z_1)$, for some $z_0 \le \a \le z_1$.
Furthermore, if $\int_{z_0}^{z_1} w > 0$, then $\a < z_1$
and $w \in B_p((\a,z_1])$.
Besides, if $w > 0$ a.e. in $[z_0,z_0']$ for some $z_0' > z_0$, then $\a = z_0$.
\end{lemma}

A similar result holds if $w$ is
a left-consistent weight in $[z_0,z_1]$.

\begin{proof}
If $\int_{z_0}^{z_1} w=0$, then $U=\emptyset=(z_1,z_1)$.
If $\int_{z_0}^{z_1} w>0$, define
$$
\a := \inf \Big\{ z \in (z_0,z_1): \, \int_{z_0}^{z} w > 0 \, \Big\} \,.
$$
It is clear that $\a< z_1$.
We show now that $U = (\a, z_1)$.
Let us fix $z\in (\a, z_1)$; then
$$
\int_{z_0}^{z} w > 0 \,,
\qquad
\Big(\int_{z_0}^{z} w \Big) \Big( \int_{z}^{z_1} w^{-1/(p-1)} \Big)^{p-1} < \infty \,,
$$
and consequently
$$
\Big( \int_{z}^{z_1} w^{-1/(p-1)} \Big)^{p-1} < \infty \,.
$$
Hence, $w \in B_p([z,z_1])$ for any $z\in (\a,z_1)$.
Then $(\a, z_1)\subset U$ and $w \in B_p((\a,z_1])$.
If $w > 0$ a.e. in $[z_0,z_0']$ for some $z_0' > z_0$, it is clear that $\a = z_0$.

If $\a=z_0$, the proof is finished.
If $\a>z_0$, then $\int_{z_0}^{\a} w = 0$ and hence
$w=0$ a.e. in $[z_0,\a]$ and
$\int_I w^{-1/(p-1)}= \infty$
for any open arc $I \subset (z_0,\a)$.
Consequently, $(z_0,\a)\cap U = \emptyset$, and $U = (\a, z_1)$.
\end{proof}

\begin{definition}
Consider $1\le p<\infty$, a compact curve $\g=[z_0,z_1]$ and a
vectorial measure $\mu=(\mu_0,\dots,\mu_{k})$ in $\g$.
We say that $\mu$ is of \emph{type} $A$ if it is finite and
strongly $p$-admissible in $\g$ and there exist points
$a_1=z_0 < a_2 < \cdots < a_{n-1} < a_n=z_1$ in $\g$
and integers $0 \le k_i^1 \le k_i^2 \le k$ $(1\le i <n)$
such that $w_j=0$ a.e. in $[a_i,a_{i+1}]$
for $k_i^2<j \le k$ if $k_i^2< k$,
and for each $1\le i <n$ we have either:

$(1)$ \ $k_i^2 = 0$,

$(2)$ \ $w_{k_i^2} \in B_p([a_i,a_{i+1}])$,

$(3)$ \ for every $k_i^1 < j\le k_i^2$,
$w_j$ is right-consistent in $[a_i,a_{i+1}]$; if $k_i^1>0$,
we also assume $a_i^+ \in \O^{(k_i^1)}$,

$(4)$ \ for every $k_i^1 < j\le k_i^2$,
$w_j$ is left-consistent in $[a_i,a_{i+1}]$; if $k_i^1>0$,
we also assume $a_{i+1}^- \in \O^{(k_i^1)}$,

$(5)$ \ for each $k_i^1 < j\le k_i^2$,
$w_j$ is either right-consistent or left-consistent in $[a_i,a_{i+1}]$; if $k_i^1>0$,
we also assume $[a_i,a_{i+1}]\subseteq \O^{(k_i^1-1)}$.
\end{definition}

\noindent
{\bf Remarks.}

{\bf 1.} The definition itemizes much more cases in order to
cover many possible behaviors of the weights. For instance, the
case $\mu=0$ in some $(a_i,a_{i+1})$ is allowed (it is contained
in the case $(1)$). In the same way, we could choose $k_i^2:=k$,
but by taking general $k_i^2$ we even allow the possibility of
different number of non-zero weights in each subarc.

{\bf 2.}
This class of measures includes many usual measures, as the Jacobi, Jacobi-Angelesco and Polacheck
weights, and measures of type $1$ or $2$ in \cite{RARP2} and \cite{APRR}.

{\bf 3.} The hypothesis $[a_i,a_{i+1}]\subseteq \O^{(k_i^1-1)}$ if
$k_i^1>0$ in $(5)$ is not very restrictive: Lemma \ref{l:8.3}
implies this condition if there exists $k_i^1 \le j_1, j_2 \le
k_i^2$ with $w_{j_1}$ right-consistent in $[a_i,a_{i+1}]$ and
$w_{j_1} > 0$ a.e. in $[a_i,a_i']$ for some $a_i' > a_i$, and
$w_{j_2}$ left-consistent in $[a_i,a_{i+1}]$ and $w_{j_2} > 0$
a.e. in $[a_{i+1}',a_{i+1}]$ for some $a_{i+1}'<a_{i+1}$ (then
$(a_i,a_{i+1}]\subseteq \O^{(j_1-1)} \subseteq \O^{(k_i^1-1)}$ and
$[a_i,a_{i+1})\subseteq \O^{(j_2-1)} \subseteq \O^{(k_i^1-1)}$).
Using the same argument, we can check that $(5)$ holds if we have
either:

$(5')$ \ for every $k_i^1 < j\le k_i^2$,
$w_j$ is right-consistent in $[a_i,a_{i+1}]$; if $k_i^1>0$,
we also assume $a_i^+ \in \O^{(k_i^1-1)}$
and $w_j >0$ a.e. in $(a_i,a_{i}')$
for some $a_i'>a_i$ and some $k_i^1 < j\le k_i^2$,

$(5'')$ \ for every $k_i^1 < j\le k_i^2$,
$w_j$ is left-consistent in $[a_i,a_{i+1}]$; if $k_i^1>0$,
we also assume $a_{i+1}^- \in \O^{(k_i^1-1)}$
and $w_j >0$ a.e. in $(a_{i+1}',a_{i+1})$
for some $a_{i+1}'<a_{i+1}$ and some $k_i^1 < j\le k_i^2$.

\medskip

\begin{teo}
\label{mult1}
Let us consider
$1\le p<\infty$ and a vectorial measure
$\mu$ of type $A$ in a compact curve $\g$.
Then the multiplication operator $\M$ is bounded in $W^{k,p}(\g,\mu)$
if and only if $\K(\g,\mu)=0$.
\end{teo}

\noindent
{\bf Remark.}
Condition $\K(\g,\mu)=0$ is easy to check in practical cases.
Propsition \ref{k=1} gives a simple characterization if $k=1$.
Although it is not possible to describe in an explicit way
the measures with $\K(\g,\mu)=0$ for any $k$
(if $\O_1\cup \cdots \cup \O_k$ is connected, then this description
is equivalent to solve the Birkhoff interpolation problem, see e.g. \cite{LJ}),
there exists a simple sufficient condition:
it is easy to check that if $\supp \mu_0$ has at least $k$ points
in each connected component of $\O_1\cup \cdots \cup \O_k$,
then $\K(\g,\mu)=0$.
Theorem \ref{mult4} below
characterizes $\K(\g,\mu)=0$ for a special kind of measures.

\begin{proof}
If $\M$ is bounded in $W^{k,p}(\g,\mu)$,
then it is well defined in $W^{k,p}(\g,\mu)$.
Since $\mu$ is a $p$-admissible vectorial measure in $\g$,
by Theorem \ref{t:8.3} we deduce that $\K(\g,\mu)=0$.

Let us assume now that $\K(\g,\mu)=0$.

We prove first that $\O^{(0)}\setminus
(\O_1\cup\cdots\cup \O_k)$ is a finite set.
Since $\O^{(0)}\subset \overline{\O_1\cup\cdots\cup \O_k}$,
by the Remark before Definition \ref{8},
it is enough to prove that
$[a_i,a_{i+1}] \cap \overline{\O_1\cup\cdots\cup \O_k}
\setminus (\O_1\cup\cdots\cup \O_k)$
is a finite set for each $1\le i <n$.

Let us fix $1\le i <n$.

If we are in case $(1)$, then
$k_i^2=0$ and $(a_i,a_{i+1}) \cap (\O_1\cup\cdots\cup \O_k)
= \emptyset$.

Consider now the case $(2)$. Then,
$w_{k_i^2} \in B_p([a_i,a_{i+1}])$; therefore $(a_i,a_{i+1})\subset \O_{k_i^2}$
and $[a_i,a_{i+1}] \cap \overline{\O_1\cup\cdots\cup \O_k} =[a_i,a_{i+1}]$;
hence $[a_i,a_{i+1}] \cap \overline{\O_1\cup\cdots\cup \O_k}
\setminus (\O_1\cup\cdots\cup \O_k) \subseteq \{a_i,a_{i+1}\}$
and it has at most two points.

Let us assume that we are in case $(5)$.
If $k_i^1 = k_i^2$, we do not require that
$w_j$ be consistent in $[a_i,a_{i+1}]$,
but then we are either in case $(1)$ if $k_i^2=0$, or $(2)$ if $k_i^2>0$:
condition $[a_i,a_{i+1}]\subseteq \O^{(k_i^2-1)}$
implies $w_{k_i^2} \in B_p([a_i,a_{i+1}])$, since
$w_j=0$ a.e. in $[a_i,a_{i+1}]$ for $k_i^2<j \le k$ if $k_i^2< k$.
Hence, without loss of generality we can assume that $k_i^1 < k_i^2$.
For each $k_i^1 < j\le k_i^2$
with $w_j$ right-consistent in $[a_i,a_{i+1}]$, by Lemma \ref{l:8.3},
$(a_i,a_{i+1}) \cap \O_j = (\a_{ij},a_{i+1})$, for some $a_i \le \a_{ij} \le a_{i+1}$.
In a similar way, for each $k_i^1 < j\le k_i^2$
with $w_j$ left-consistent in $[a_i,a_{i+1}]$, we have that
$(a_i,a_{i+1}) \cap \O_j = (a_i,\b_{ij})$, for some $a_i \le \b_{ij} \le a_{i+1}$.
If we define
$$
\begin{aligned}
\a_i & :=\min \big\{\a_{ij}: \, k_i^1 < j\le k_i^2
\text{ and $w_j$ is right-consistent in $[a_i,a_{i+1}]$} \, \big\}\,,
\\
\b_i & :=\max \big\{\b_{ij}: \, k_i^1 < j\le k_i^2
\text{ and $w_j$ is left-consistent in $[a_i,a_{i+1}]$} \, \big\}\,,
\end{aligned}
$$
then $(a_i,a_{i+1}) \cap (\O_{k_i^1+1} \cup\cdots\cup \O_k)
=(a_i,\b_i) \cup (\a_i,a_{i+1})$.

If $\b_i \ge \a_i$, then
$(a_i,a_{i+1}) \setminus \{\a_i\} \subseteq
(a_i,a_{i+1}) \cap (\O_{1} \cup\cdots\cup \O_k)$,
and consequently
$[a_i,a_{i+1}] \cap \overline{\O_1\cup\cdots\cup \O_k}
\setminus (\O_1\cup\cdots\cup \O_k) \subseteq \{a_i, \a_i, a_{i+1}\}$.

If $\b_i < \a_i$ and $k_i^1>0$, then
$(\b_i , \a_i) \cap (\O_{k_i^1+1}\cup\cdots\cup \O_k) = \emptyset$ and
$[a_i,a_{i+1}] \subseteq \O^{(k_i^1-1)}
\subseteq \overline{\O_{k_i^1}\cup\cdots\cup \O_k}$,
and hence
$(\b_i , \a_i) \subseteq \O_{k_i^1}$.
Therefore,
$(a_i,a_{i+1}) \setminus \{\b_i , \a_i\} \subseteq \O_{k_i^1}\cup\cdots\cup \O_k$,
and
$[a_i,a_{i+1}] \cap \overline{\O_1\cup\cdots\cup \O_k}
\setminus (\O_1\cup\cdots\cup \O_k) \subseteq \{a_i, \b_i, \a_i, a_{i+1}\}$.

If $\b_i < \a_i$ and $k_i^1=0$,
then $w_j=0$ a.e. in $(\b_i , \a_i)$ for every $1\le j \le k$,
$(\b_i , \a_i) \cap (\O_{1} \cup\cdots\cup \O_k) = \emptyset$
and $(\b_i , \a_i) \cap \overline{\O_1\cup\cdots\cup \O_k} = \emptyset$.
Consequently,
$[a_i,a_{i+1}] \cap \overline{\O_1\cup\cdots\cup \O_k}
\setminus (\O_1\cup\cdots\cup \O_k) \subseteq \{a_i, \b_i, \a_i, a_{i+1}\}$.

Consider now the case $(3)$.
Then
$(a_i,a_{i+1}) \cap (\O_{k_i^1+1} \cup\cdots\cup \O_k) = (\a_i,a_{i+1})$.

If $k_i^1>0$, then $a_i^+ \in \O^{(k_i^1)}$.
By the Remark before Definition \ref{8},
there exists $a_i' > a_i$ such that
$(a_i,a_{i}') \subset \O_{k_i^1+1} \cup\cdots\cup \O_k$,
and consequently, $\a_i=a_i$.
Then,
$(a_i,a_{i+1}) \subseteq \O_{k_i^1+1}\cup\cdots\cup \O_k$,
and
$[a_i,a_{i+1}] \cap \overline{\O_1\cup\cdots\cup \O_k}
\setminus (\O_1\cup\cdots\cup \O_k) \subseteq \{a_i, a_{i+1}\}$.

If $k_i^1=0$ and $\a_i=a_i$, we also have this inclusion.

If $k_i^1=0$ and $\a_i > a_i$,
then $w_j=0$ a.e. in $(a_i , \a_i)$ for every $1\le j \le k$,
$(a_i , \a_i) \cap (\O_{1} \cup\cdots\cup \O_k) = \emptyset$
and $(a_i , \a_i) \cap \overline{\O_1\cup\cdots\cup \O_k} = \emptyset$.
Consequently,
$[a_i,a_{i+1}] \cap \overline{\O_1\cup\cdots\cup \O_k}
\setminus (\O_1\cup\cdots\cup \O_k) \subseteq \{\a_i, a_{i+1}\}$.

The case $(4)$ are similar to $(3)$.

\smallskip

Since
$\O^{(0)}\setminus (\O_1\cup\cdots\cup \O_k)$ is a finite set
and $\K(\g,\mu)=0$,
we conclude that $(\g,\mu)\in \C_0$
(see Remark 1 after Definition \ref{d:claseC}).

By Lemma \ref{l:8.2}, we just need to show that
there exists a positive constant $c$ such that
\begin{equation}
\label{eq:1}
\|f^{(j-1)}\|_{L^{p}([a_i,a_{i+1}],\mu_j)}
\le c\, \|f\|_{W^{k,p}(\g,\mu)} \,,
\end{equation}
for every $1\le i <n$, $1\le j\le k_i ^2$ (if $k_i ^2 \ge 1$) and
$f\in V^{k,p}(\g,\mu)$.

Let us fix $1\le i<n$.

If we are in case $(1)$, then $k_i^2=0$ and there is nothing to prove.
We can assume now that $k_i ^2 \ge 1$.

Consider now the case $(2)$. Then,
$w_{k_i^2} \in B_p([a_i,a_{i+1}])$; therefore
$[a_i,a_{i+1}] \subseteq \O^{(k_i^2 -1)}$
and consequently
$[a_i,a_{i+1}] \subseteq \O^{(j -1)}$ for $1\le j\le k_i ^2$;
by Theorem \ref{t:4.1} we have directly
$$
\|f^{(j-1)}\|_{L^{p}([a_i,a_{i+1}],\mu_j)}
\le c_1\, \|f^{(j-1)}\|_{L^{\infty}([a_i,a_{i+1}])}
\le c_2\, \|f\|_{W^{k,p}(\g,\mu)} \,,
$$
for every $f\in V^{k,p}(\g,\mu)$ and $1\le j\le k_i^2$,
since $(\g,\mu)\in\C_0$.

Let us assume now that we are in case $(5)$.
As above, without loss of generality we can assume that $k_i^1 < k_i^2$.

If $w_j=0$ a.e. in $[a_i,a_{i+1}]$ for some $k_i^1 < j\le k_i^2$, then we have
$$
\|f^{(j-1)}\|_{L^{p}([a_i,a_{i+1}],w_j)}
\le c\,\|f\|_{W^{k,p}(\g,\mu)}\,,
$$
for every positive constant $c$.

Fix $k_i^1 < j\le k_i^2$ with
$w_j$ right-consistent in $[a_i,a_{i+1}]$
and $\int_{a_{i}}^{a_{i+1}} w_{j} > 0$.
Muckenhoupt inequality gives that
$$
\begin{aligned}
c_3\,\|f^{(j-1)}(a_{i+1}) - f^{(j-1)}\|_{L^{p}([a_i,a_{i+1}],w_j)}
& \le
\|f^{(j)}\|_{L^{p}([a_i,a_{i+1}],w_j)}\,,
\\
c_4\,\|f^{(j-1)}\|_{L^{p}([a_i,a_{i+1}],w_j)}
& \le
\|f^{(j)}\|_{L^{p}([a_i,a_{i+1}],w_j)}+|f^{(j-1)}(a_{i+1})| \,.
\end{aligned}
$$
Since $\int_{a_{i}}^{a_{i+1}} w_{j} > 0$,
we deduce that $w_{j} \in B_p([a_{i+1}', a_{i+1}])$
for some $a_{i+1}'< a_{i+1}$,
by Lemma \ref{l:8.3}.
Hence $a_{i+1}$ is $(j-1)$-left-regular
(see Remark 1(a) after Definition \ref{d:regular}).
Since $(\g,\mu)\in\C_0$, Theorem \ref{t:4.1} gives
$$
|f^{(j-1)}(a_{i+1})|
\le c_5\,\|f\|_{W^{k,p}(\g,\mu)}\,,
$$
and we conclude
$$
\|f^{(j-1)}\|_{L^{p}([a_i,a_{i+1}],w_j)}
\le c_6\,\|f\|_{W^{k,p}(\g,\mu)}\,.
$$
If we fix $k_i^1 < j\le k_i^2$ with
$w_j$ left-consistent in $[a_i,a_{i+1}]$
and $\int_{a_{i}}^{a_{i+1}} w_{j} > 0$,
we obtain a similar inequality.
Consequently,
$$
\|f^{(j-1)}\|_{L^{p}([a_i,a_{i+1}],w_j)}
\le c_7\,\|f\|_{W^{k,p}(\g,\mu)}\,,
$$
for every $k_i^1< j\le k_i^2$.

For each $k_i^1< j\le k_i^2$, we define
$$
\begin{aligned}
\a_i^j & :=\min \big\{\a_{im}: \, j < m\le k_i^2
\text{ and $w_m$ is right-consistent in $[a_i,a_{i+1}]$} \, \big\}\,,
\\
\b_i^j & :=\max \big\{\b_{im}: \, j < m\le k_i^2
\text{ and $w_m$ is left-consistent in $[a_i,a_{i+1}]$} \, \big\}\,.
\end{aligned}
$$
Then
$(a_i,a_{i+1}) \cap (\O_{j+1}\cup\cdots\cup \O_k)=(a_i,\b_i^j) \cup (\a_i^j,a_{i+1})$
and
$[a_i,a_{i+1}] \cap \O^{(j)} \subseteq
[a_i,a_{i+1}] \cap \overline{\O_{j+1}\cup\cdots\cup \O_k} =[a_i,\b_i] \cup [\a_i,a_{i+1}]$;
consequently $[a_i,a_{i+1}] \cap \O^{(j)}$ has at most two connected components.

Let us assume that $[a_i,a_{i+1}] \cap \O^{(j)}$ has two connected components $\L_1, \L_2$.
Since $\mu$ is a strongly $p$-admissible vectorial measure,
we have that
$\supp (\mu_j|_{\L_l})_s \subseteq \supp (\mu_j^*|_{\L_l})
\subseteq \L_l \subseteq \O^{(j)} \subseteq \O^{(j-1)}$ ($l=1,2$).
The smaller arc $K_j^l$ which contains to $\supp (\mu_j|_{\L_l})_s$, is
a compact arc contained in $\L_l \subseteq \O^{(j-1)}$ ($l=1,2$).
Hence, Theorem \ref{t:4.1} gives
$$
\|f^{(j-1)}\|_{L^{p}([a_i,a_{i+1}],(\mu_j)_s)}
\le c_8\,\|f^{(j-1)}\|_{L^{\infty}(K_j^1 \cup K_j^2)}
\le c_9\,\|f\|_{W^{k,p}(\g,\mu)}\,,
$$
for $k_i^1< j\le k_i^2$, since $(\g,\mu)\in\C_0$.
If $[a_i,a_{i+1}] \cap \O^{(j)}$ is connected,
we obtain a similar inequality.

Consequently,
\begin{equation}
\label{eq:1.1}
\|f^{(j-1)}\|_{L^{p}([a_i,a_{i+1}],\mu_j)}
\le c_{10}\,\|f\|_{W^{k,p}(\g,\mu)}\,,
\end{equation}
for $k_i^1< j\le k_i^2$.
If $k_i^1=0$, we have finished the proof.

Assume now that $k_i^1>0$ and fix $1\le j\le k_i^1$.
Since $[a_i,a_{i+1}] \subseteq \O^{(k_i^1 -1)} \subseteq \O^{(j -1)}$
and $(\g,\mu)\in\C_0$,
Theorem \ref{t:4.1} gives again
$$
\|f^{(j-1)}\|_{L^{p}([a_i,a_{i+1}],\mu_j)} \le
c_{11}\,\|f^{(j-1)}\|_{L^{\infty}([a_i,a_{i+1}])}
\le c_{12}\,\|f\|_{W^{k,p}(\g,\mu)} \,,
$$
for all $f\in V^{k,p}(\g,\mu)$ and $1\le j\le k_i^1$.

Therefore $(\ref{eq:1})$ holds
for every $1\le j\le k_i^2$ and $f\in V^{k,p}(\g,\mu)$.

\smallskip

We consider now the case $(3)$.
As above, without loss of generality we can assume that $k_i^1 < k_i^2$.
The argument to obtain
$(\ref{eq:1.1})$ is similar to the one in case $(5)$;
here, $[a_i,a_{i+1}] \cap \O^{(j)}$ is connected
for $k_i^1< j\le k_i^2$.
If $k_i^1=0$, we have finished the proof.

Assume now that $k_i^1>0$.
We have seen that
$(a_i,a_{i+1}) \cap (\O_{k_i^1+1} \cup\cdots\cup \O_k) = (\a_i,a_{i+1})$
in case $(3)$.
Since $a_i^+ \in \O^{(k_i^1)}$,
by the Remark before Definition \ref{8},
there exists $a_i' > a_i$ such that
$(a_i,a_{i}') \subset \O_{k_i^1+1} \cup\cdots\cup \O_k$,
and consequently, $\a_i=a_i$.
Then,
$(a_i,a_{i+1}) \subseteq \O_{k_i^1+1} \cup\cdots\cup \O_k
\subseteq \O^{(k_i^1)}$.

Without loss of generality we can assume that
$$
k_i^2 =\min \big\{m \le k : \, w_j=0
\text{ a.e. in $[a_i,a_{i+1}]$ for } m < j \le k \, \big\}\,;
$$
then $\int_{a_{i}}^{a_{i+1}} w_{k_i^2} > 0$.

Since $\int_{a_{i}}^{a_{i+1}} w_{k_i^2} > 0$,
we deduce that $w_{k_i^2} \in B_p((\a_{ik_i^2} , a_{i+1}])$
for some $\a_{ik_i^2} < a_{i+1}$,
by Lemma \ref{l:8.3}.
Then, $a_{i+1}^- \in \O^{(k_i^2-1)} \subseteq \O^{(k_i^1)}$.
Recall that $a_i^+ \in \O^{(k_i^1)}$.
Hence, $[a_i,a_{i+1}] \subseteq \O^{(k_i^1)}$.

Fix $1\le j\le k_i^1$.
Since $[a_i,a_{i+1}] \subseteq \O^{(k_i^1)} \subseteq \O^{(j -1)}$
and $(\g,\mu)\in\C_0$,
Theorem \ref{t:4.1} gives again
$$
\|f^{(j-1)}\|_{L^{p}([a_i,a_{i+1}],\mu_j)} \le
c_{11}\,\|f^{(j-1)}\|_{L^{\infty}([a_i,a_{i+1}])}
\le c_{12}\,\|f\|_{W^{k,p}(\g,\mu)} \,,
$$
for all $f\in V^{k,p}(\g,\mu)$ and $1\le j\le k_i^1$.

Therefore $(\ref{eq:1})$ holds
for every $1\le j\le k_i ^2$ and $f\in V^{k,p}(\g,\mu)$.

\smallskip

The case $(4)$ is similar to $(3)$.

\smallskip

Then $(\ref{eq:1})$ holds
for every $1\le i <n$, $1\le j\le k_i ^2$ (if $k_i^2 \ge 1$) and $f\in V^{k,p}(\g,\mu)$.
Consequently, Lemma \ref{l:8.2} finishes the proof.
\end{proof}

As we said, if $k=1$, there exists a simple characterization of $\K(\g,\mu)=0$.
This result is interesting since in the applications of Sobolev orthogonal polynomials
usually appears the case $k=1$.

\begin{prop}
\label{k=1}
Let us  consider $1\le p < \infty$ and a $p$-admissible vectorial measure
$\mu=(\mu_0,\mu_1)$ in a curve $\g$.
Then $\K(\g,\mu)=0$ if and only $\mu_0(A)>0$
for every connected component $A$ of $\O^{(0)}$.
\end{prop}

\begin{proof}
If there exists a connected component $A$ of $\O^{(0)}$
with $\mu_0(A)=0$, then $\chi_{_{A}}\in \K(\g,\mu)$
and $\K(\g,\mu)\neq 0$.

Assume now that $\mu_0(A)>0$
for every connected component $A$ of $\O^{(0)}$.
Fix any connected component $A$ of $\O^{(0)}$.
Note that every interior point of $A$ (as a subset of $\g$)
belongs to $\O_1$.
Consequently, if $f\in \K(\g,\mu)$, then $f'=0$ a.e. in the interior of $A$,
and $f$ is constant in $A$.
Since $\mu_0(A)>0$, then $f=0$ in $A$.
Hence, $f=0$ in $\O^{(0)}$ and $\K(\g,\mu)=0$.
\end{proof}

It could be interesting to check Theorem \ref{mult1}
in some particular case.
Consider again the example after Definition \ref{heuristic}.
Since $\K (\g,\mu)=\span \{x, (x_+)^2\}$,
$\M$ is not bounded, but in order to make $\M$ bounded in the Sobolev space,
it is sufficient to replace
$\mu_0=\d_0$ by a measure with two more Dirac's deltas.

But, what is the minimum amount of deltas that we need in order
to have the multiplication operator bounded in the closure of the
space of polynomials $\PP^{3,p}([-1,1],\mu)$ with the Sobolev norm?

We need at least one:
In order to have a norm in the space of polynomials
(which is equivalent to the existence of Sobolev orthogonal polynomials),
we need to replace
$\mu_0=\d_0$ by another measure with one more delta
(the polynomials in $\K([-1,1],\mu)$ are just the span of $x$).

But, in fact, we need two, since, as usual,
the Sobolev space is the closure of
the space of polynomials $\PP^{3,p}([-1,1],\mu)$ with the Sobolev norm.

Therefore our definition of Sobolev space gives the sharp result:
we see that even if $\|\cdot \|_{W^{k,p}(\g,\mu)}$ is a norm in $\PP$,
$\M$ can be not bounded in $\PP^{k,p}(\g,\mu)$;
we need $\|\cdot \|_{W^{k,p}(\g,\mu)}$
to be a norm in $W^{k,p}(\g,\mu)$.

\smallskip

After this general result, we can deduce three practical consequences.

\begin{definition}
A function $u$ in a compact curve $\g=[z_0,z_1]$ is
\emph{piecewise monotone} if there exist points
$b_1=z_0 < b_2 < \cdots < b_{m-1} < b_m=z_1$ in $\g$
such that $u$ is a monotone function in $[b_i,b_{i+1}]$
for each $1\le i <m$.
\end{definition}

\begin{definition}
\label{d:typeB}
Consider $1\le p<\infty$ and a
vectorial measure $\mu=(\mu_0,\dots,\mu_{k})$ in a compact curve $\g$.
We say that $\mu$ is of \emph{type} $B$ if it is finite and
strongly $p$-admissible in $\g$ and
$w_j$ is comparable to a piecewise monotone function
for any $1 \le j \le k$.
\end{definition}

\noindent {\bf Remarks.}

{\bf 1.}
By monotone we mean non-strictly monotone;
hence, it is possible to have $w_j=0$ in some arc.

{\bf 2.}
The partition in arcs can be different for each $w_j$.

\smallskip

\begin{teo}
\label{mult2}
Let us consider
$1\le p<\infty$ and a vectorial measure
$\mu$ of type $B$ in a compact curve $\g$.
Then the multiplication operator $\M$ is bounded in $W^{k,p}(\g,\mu)$
if and only if $\K(\g,\mu)=0$.
\end{teo}

\begin{proof}
We just need to prove that $\mu$ is a measure of type $A$ in $\g$,
with cases $(1)$ and/or $(5)$;
hence, Theorem \ref{mult1} finishes the proof.

For each $1 \le j \le k$, there exist points
$b_1^j=z_0 < b_2^j < \cdots < b_{m^j-1}^j < b_{m^j}^j =z_1$ in $\g$
such that $w_j$ is comparable to a monotone function in $[b_i^j,b_{i+1}^j]$
for each $1\le i <m^j$.
Splitting in two subarcs some arcs if it is necessary,
without loss of generality we can require also that in each arc $[b_i^j,b_{i+1}^j]$,
we have either $w_j=0$ a.e. or $w_j>0$ a.e.

We consider the points
$a_1=z_0 < a_2 < \cdots < a_{n-1} < a_n=z_1$ in $\g$, which are
the ordered points in the set $\{b_i^j\}_{1\le i <m^j, \, 1 \le j \le k}$.
Consequently,
for any fixed $1\le i <n$ and $1 \le j \le k$,
$w_j$ is comparable to a monotone function in $[a_i,a_{i+1}]$
and we have either $w_j=0$ a.e. or $w_j>0$ a.e. in $[a_i,a_{i+1}]$.

For any $1\le i< n$, we define $k_i^1 := 0$ and
$$
k_i^2 :=\min \big\{m \le k : \, w_j=0
\text{ a.e. in $[a_i,a_{i+1}]$ for } m < j \le k \, \big\}\,.
$$
If $k_i^2= 0$, then the arc $[a_i,a_{i+1}]$ is in the case $(1)$
in the definition of measure of type $A$.
If $k_i^2> 0$, then the arc $[a_i,a_{i+1}]$ is in the case $(5)$:
the Remark after Definition \ref{d:comparable}
gives that if $w_j$ is comparable to a
non-decreasing function in $[a_i,a_{i+1}]$
then it is right-consistent, and
if $w_j$ is comparable to a
non-increasing function in $[a_i,a_{i+1}]$
then it is left-consistent.
If $w_j=0$ a.e. in $[a_i,a_{i+1}]$,
then it is both right and left-consistent.

Then $\mu$ is a measure of type $A$ in $\g$
and Theorem \ref{mult1} gives the result.
\end{proof}

Theorem \ref{mult2} gives the following direct result.

\begin{teo}
\label{mult3}
Let us consider
$1\le p<\infty$ and a finite vectorial measure
$\mu$ in a compact curve $\g$ such that
$d\mu_j=w_j ds$ and
$w_j$ is comparable to a piecewise monotone function
for any $1 \le j \le k$.
Then the multiplication operator $\M$ is bounded in $W^{k,p}(\g,\mu)$
if and only if $\K(\g,\mu)=0$.
\end{teo}

Theorem \ref{mult3} is an improvement of \cite[Theorem 4.3]{R4},
since here we consider Sobolev spaces in curves instead of intervals;
furthermore, it is also an improvement of \cite[Theorem 4.3]{R4} in the case of
intervals: In \cite{R4} we consider a different Sobolev space (in an interval $I$),
which we denote by $W_{ko}^{k,p}(I,\mu)$, verifying
$\PP^{k,p}(I,\mu) \subseteq W^{k,p}(I,\mu) \subseteq W_{ko}^{k,p}(I,\mu)$.
Since we usually have
$\PP^{k,p}(I,\mu) = W^{k,p}(I,\mu) \neq W_{ko}^{k,p}(I,\mu)$
(see \cite[Theorem 6.1]{APRR}),
it is obvious that it is better to work with $W^{k,p}(I,\mu)$ in order
to obtain results about the multiplication operator in $\PP^{k,p}(I,\mu)$.
(The advantage of $W_{ko}^{k,p}(I,\mu)$ is that it can be defined in a simpler and
faster way than $W^{k,p}(I,\mu)$.)

\begin{proof}
We prove that, with our hypothesis, $\mu$ is strongly $p$-admissible in $\g$;
in fact, we prove that $\mu_j^*=0$ for every $1 \le j \le k$.
Let us fix $1 \le j \le k$.
Since $d\mu_j^* := d\mu_j - w_j \chi_{_{\Omega_{j}}} \! ds$
and $\mu_j$ is absolutely continuous in $\g$,
$d\mu_j^* := w_j (1 - \chi_{_{\Omega_{j}}} ) ds$, and
we just need to prove $w_j=0$ a.e. in $\g \setminus \Omega_{j}$.
This fact is a consequence of Lemma \ref{l:8.3}:
using the notation in the proof of Theorem \ref{mult2},
since for each $1\le i< n$, we have that $w_j$ is either
right or left-consistent in $[a_i,a_{i+1}]$,
and we have either $w_j=0$ a.e. or $w_j>0$ a.e. in $[a_i,a_{i+1}]$,
Lemma \ref{l:8.3} gives either $w_j=0$ a.e. or $w_j\in B_p((a_i,a_{i+1}))$.
Consequently,
for each $1\le i< n$, we have either $w_j=0$ a.e. in $[a_i,a_{i+1}]$
or $(a_i,a_{i+1}) \subseteq \O_j$;
in both cases $w_j=0$ a.e. in $(\g \setminus \O_j) \cap [a_i,a_{i+1}]$,
and consequently $w_j=0$ a.e. in $\g \setminus \O_j$.
Then, $\mu_j^*=0$ for every $1 \le j \le k$
and hence $\mu$ is strongly $p$-admissible in $\g$.
Therefore, $\mu$ is of type B, and Theorem \ref{mult2} finishes the proof.
\end{proof}

It is usual that the behavior of a weight is ``similar" to some power, in some sense;
the following definition deals with this case.

\begin{definition}
\label{d:typeC}
Consider $1\le p<\infty$, a compact curve $\g=[a_1,a_4]$ with
$|z-a_1| \asymp |\g^{-1} (z) -\g^{-1}(a_1)|$ for a.e. $z \in [a_1,a_2]$ and
$|z-a_4| \asymp |\g^{-1} (z) -\g^{-1}(a_4)|$ for a.e. $z \in [a_3,a_4]$,
where $\g^{-1}$ denotes the inverse of some
parametrization with $\g' \in L^{\infty}([a_1,a_4])$,
and a vectorial measure $\mu=(\mu_0,\dots,\mu_{k})$ in $\g$.
We say that $\mu$ is of \emph{type} $C$ if it is finite and
strongly $p$-admissible in $\g$, and we have:

$(1)$ \ $w_{k} \in B_p((a_1,a_4))$,

$(2)$ \ for each $1\le j \le k$ we have either:

$\quad (2.1)$ \ $w_j$ is comparable to a monotone function in $[a_1,a_2]$,

$\quad (2.2)$ \ $w_j \in B_p([a_1,a_2])$,

$\quad (2.3)$ \ $c_1 |z-a_1|^{p-1} \le w_j(z) \le c_2 |z-a_1|^{\b_j^0}$ a.e. in
$[a_1,a_2]$, with $\b_j^0 > -1$,

$\quad (2.4)$ \ $c_1 |z-a_1|^{\a_j^0} \le w_j(z) \le c_2 |z-a_1|^{\a_j^0-p}$ a.e. in
$[a_1,a_2]$, with $\a_j^0 > p-1$,

$(3)$ \ for each $1\le j \le k$ we have either:

$\quad (3.1)$ \ $w_j$ is comparable to a monotone function in $[a_3,a_4]$,

$\quad (3.2)$ \ $w_j \in B_p([a_3,a_4])$,

$\quad (3.3)$ \ $c_1 |z-a_4|^{p-1} \le w_j(z) \le c_2 |z-a_4|^{\b_j^1}$ a.e.
in $[a_3,a_4]$, with $\b_j^1 > -1$,

$\quad (3.4)$ \ $c_1 |z-a_4|^{\a_j^1} \le w_j(z) \le c_2 |z-a_4|^{\a_j^1-p}$ a.e.
in $[a_3,a_4]$, with $\a_j^1 > p-1$.
\end{definition}

\noindent
{\bf Remark.}
In $(2.3)$ and $(2.4)$ we require $c_1 |z-a_1|^{\a} \le w_j(z)$
a.e. in $[a_1,a_2]$, with $\a \ge p-1$.
We also cover the other cases in $(2.2)$:
If $c_1 |z-a_1|^{\a} \le w_j(z)$ a.e. in
$[a_1,a_2]$, with $\a < p-1$, then $w_j \in B_p([a_1,a_2])$
(without any upper bound for $w_j$).

The same holds in $[a_3,a_4]$.

\begin{teo}
\label{mult4}
Let us consider
$1\le p<\infty$ and a vectorial measure
$\mu$ of type $C$ in a compact curve $\g$.
Then the following conditions are equivalent:

$(i)$ The multiplication operator $\M$ is bounded in $W^{k,p}(\g,\mu)$.

$(ii)$ $\K(\g,\mu)=0$.

Furthermore, if $\int_\g w_1>0$, conditions $(i)$ and $(ii)$ are also equivalent to

$(iii)$ $\mu_0(\g)>0$.
\end{teo}

\begin{proof}
In order to prove the equivalence between $(i)$ and $(ii)$,
we just need to prove that $\mu$ is a measure of type $A$ in $\g$
and to apply Theorem \ref{mult1}.
We deal with the case $1< p<\infty$ in order to write $1/(p-1)$ as a real number;
the case $p=1$ is similar.

Let us define
$k_1^1 = k_2^1 = k_3^1 = 0$ and $k_1^2 = k_2^2 = k_3^2 = k$.
Without loss of generality we can assume that $a_2 < a_3$;
then $a_1 < a_2 < a_3 < a_4$.

Since $w_{k} \in B_p([a_2,a_{3}])$, this arc verifies hypothesis
$(2)$ of measures of type $A$.

We prove now that the arc $[a_1,a_{2}]$ verifies hypothesis
$(5)$ of measures of type $A$. Let us fix $0 < j \le k$.

In case $(2.1)$, $w_j$ is comparable to a monotone function in $[a_1,a_2]$, and
then $w_j$ is right-consistent or left-consistent in $[a_1,a_{2}]$.

In case $(2.2)$, $w_j$ is also right-consistent in $[a_1,a_2]$
(and left-consistent in this interval), since
$$
\L_{p,[a_1,a_{2}]}^+ (w_j, w_j)
= \sup_{a_1<z<a_2} \Big( \int_{a_1}^z w_j \Big)
\Big(\int_{z}^{a_2} w_j^{-1/(p-1)} \Big)^{p-1}
\le
\Big( \int_{a_1}^{a_2} w_j \Big)
\Big(\int_{a_1}^{a_2} w_j^{-1/(p-1)} \Big)^{p-1} < \infty \,.
$$

In case $(2.4)$,
let us consider the arc-length parametrization
$\g: [0,l] \longrightarrow [a_1,a_{2}]$. We have
$$
\begin{aligned}
\int_{a_1}^z w_j(\z) |d\z| \, &
\Big(\int_{z}^{a_2} w_j(\z)^{-1/(p-1)} |d\z| \Big)^{p-1}
\le
\int_{a_1}^{z} c_2 |\z-a_1|^{\a_j^0-p} |d\z| \,
\Big(\int_{z}^{a_2} c_1^{-1/(p-1)}  |\z-a_1|^{-\a_j^0/(p-1)} |d\z| \Big)^{p-1}
\\
& \le
c \int_{a_1}^{z} |\g^{-1}(\z)-\g^{-1}(a_1)|^{\a_j^0-p} |d\z| \,
\Big(\int_{z}^{a_2} |\g^{-1}(\z)-\g^{-1}(a_1)|^{-\a_j^0/(p-1)} |d\z| \Big)^{p-1}
\\
& =
c \int_{0}^{t} s^{\a_j^0-p} |\g'(s)| \, ds \,
\Big(\int_{t}^{l} s^{-\a_j^0/(p-1)} |\g'(s)| \, ds \Big)^{p-1}
\\
& \le
c' \int_{0}^{t} s^{\a_j^0-p} ds \,
\Big(\int_{t}^{l} s^{-\a_j^0/(p-1)} ds \Big)^{p-1}
\\
& =
c' \,\, \frac{t^{\a_j^0-p+1}}{\a_j^0-p+1} \,
\Big(\frac{p-1}{\a_j^0-p+1} \big(t^{(-\a_j^0+p-1)/(p-1)} - l^{(-\a_j^0+p-1)/(p-1)}\big)
\Big)^{p-1}
\\
& =
c''
\Big(1 - (t/l)^{(\a_j^0-p+1)/(p-1)} \Big)^{p-1}
\le c'' ,
\end{aligned}
$$
for any $z\in (a_1,a_2)$ (i.e. $t \in (0,l)$),
since $\a_j^0 > p-1$.
Hence,
$\L_{p,[a_1,a_{2}]}^+ (w_j, w_j) < \infty$ and
$w_j$ is right-consistent in $[a_1,a_{2}]$.

Case $(2.3)$ is similar to $(2.4)$.

Therefore, the arc $[a_1,a_{2}]$ verifies hypothesis
$(5)$ of measures of type $A$.

A similar argument proves that the arc $[a_3,a_{4}]$ also verifies hypothesis
$(5)$ of measures of type $A$.
Then, $\mu$ is a measure of type $A$ in $\g$,
and $(i)$ is equivalent to $(ii)$.

\smallskip

We prove now the equivalence of $(ii)$ and $(iii)$:

\smallskip

If $\mu_0(\g)=0$, then
$$
\big\| 1 \big\|_{W^{k,p}(\g,\mu)}^p
= \int_\g | 1 |^p d\mu_0 = 0 \,,
$$
and consequently $1 \in \K(\g,\mu) \neq 0$.

\smallskip

If $\int_\g w_1 >0$ and $\mu_0(\g)>0$, let us consider $f \in \K(\g,\mu)$.
Condition $(1)$ gives $\O_k=(a_1,a_4)$, and hence
Proposition \ref{p:fj} implies $f \in \PP_{k-1}$.
Since $\int_\g w_1 > 0$ and
$\int_\g | f' |^p w_1 = 0$, we obtain $f' = 0$ in infinitely many points in $\g$;
consequently, $f$ is constant.
Since $\mu_0(\g)>0$ and $| f |^p \mu_0(\g) = \int_\g | f |^p d\mu_0 = 0$,
we have that $f = 0$ in $\g$.
Then $\K(\g,\mu) = 0$.
\end{proof}


\begin{thebibliography}{99}

\bibitem{APRR} Alvarez, V., Pestana, D.,
Rodr\'{\i}guez, J. M., Romera, E.: Weighted Sobolev spaces on
curves. J. Approx. Theory 119, 41-85 (2002)

\bibitem{BFM} Branquinho, A., Foulqui\'e, A., Marcell\'an, F.:
Asymptotic behavior of Sobolev type orthogonal polynomials on a
rectifiable Jordan curve or arc. Constr. Approx. 18, 161-182 (2002)

\bibitem{CM} Cachafeiro, A., Marcell\'an, F.:
Orthogonal polynomials of Sobolev type on the unit circle.
J. Approx. Theory 78, 127-146 (1994)

\bibitem{CD} Castro, M., Dur\'an, A. J.:
Boundedness properties for Sobolev inner products.
J. Approx. Theory 122, 97-111 (2003)

\bibitem{EL} Everitt, W. N., Littlejohn, L. L.: The density of polynomials in
a weighted Sobolev space. Rendiconti di Matematica, Serie VII,
10, 835-852 (1990)

\bibitem{ELW1} Everitt, W. N., Littlejohn, L. L., Williams, S. C.:
Orthogonal polynomials in weighted Sobolev spaces. In:
Lecture Notes in Pure and Applied Mathematics, 117,
Marcel Dekker, 53-72 (1989)

\bibitem{ELW2} Everitt, W. N., Littlejohn, L. L., Williams, S. C.:
Orthogonal polynomials and approximation in Sobolev spaces.
J. Comput. Appl. Math. 48, 69-90 (1993)

\bibitem{FMP} Foulqui\'e, A., Marcell\'an, F., Pan, K.:
Asymptotic behavior of Sobolev-type orthogonal polynomials on the unit
circle. J. Approx. Theory 100, 345-363 (1999)

\bibitem{IKNS1} Iserles, A., Koch, P. E., Norsett, S. P., Sanz-Serna, J. M.:
Orthogonality and approximation in a Sobolev space. In: Algorithms for
Approximation. J. C. Mason and M. G. Cox, Chapman \& Hall, London (1990)

\bibitem{IKNS2} Iserles, A., Koch, P. E., Norsett, S. P., Sanz-Serna, J. M.:
On polynomials orthogonal with respect to certain Sobolev inner
products. J. Approx. Theory 65, 151-175 (1991)

\bibitem{Ku} Kufner, A.: Weighted Sobolev Spaces. Teubner
Verlagsgesellschaft. Teubner-Texte zur Mathematik (Band 31), Leipzig (1980).
Also published  by John Wiley \& Sons, New York (1985)

\bibitem{KO} Kufner, A., Opic, B.: How to define reasonably
Weighted Sobolev Spaces. Comment. Math. Univ. Carolinae 25(3), 537-554 (1984)

\bibitem{LP} L\'opez Lagomasino, G., Pijeira, H.: Zero location and
$n$-th root asymptotics of Sobolev orthogonal polynomials. J.
Approx. Theory 99, 30-43 (1999)

\bibitem{LPP} L\'opez Lagomasino, G., Pijeira, H., P\'erez, I:
Sobolev orthogonal polynomials in the complex plane. J. Comp.
Appl. Math. 127, 219-230 (2001)

\bibitem{LJ} Lorentz, G. G., Jetter, K.: Riemenschneider, Birkhoff interpolation.
Addison-Wesley Pub. (1984)

\bibitem{M-F} Mart{\'\i}nez-Finkelshtein, A.:
Bernstein-Szeg\"o's theorem for Sobolev orthogonal polynomials.
Constr. Approx. 16, 73-84 (2000)

\bibitem{M} Maz'ja, V. G.: Sobolev spaces. Springer-Verlag, New York (1985)

\bibitem{Mu2} Muckenhoupt, B.: Hardy's inequality with weights.
Studia Math. 44, 31-38 (1972)

\bibitem{PQRT1}
Portilla, A., Quintana, Y., Rodr{\'\i}guez, J. M., Tour{\'\i}s, E.:
Weierstrass' Theorem with weights.
J. Approx. Theory 127, 83-107 (2004)

\bibitem{PQRT2}
Portilla, A., Quintana, Y., Rodr{\'\i}guez, J. M., Tour{\'\i}s, E.:
Weighted Weierstrass' Theorem with first derivatives.
J. Math. Anal. Appl. 334, 1167-1198 (2007)

\bibitem{PQRT3}
Portilla, A., Quintana, Y., Rodr{\'\i}guez, J. M., Tour{\'\i}s, E.:
Weierstrass' Theorem in weighted Sobolev spaces with $k$ derivatives.
Rocky Mount. J. Math. 37, 1989-2024 (2007)

\bibitem{RARP1} Rodr\'{\i}guez, J. M., Alvarez, V., Romera, E., Pestana, D.:
Generalized weighted Sobolev spaces and applications to Sobolev
orthogonal polynomials I.
Acta Appl. Math. 80, 273-308 (2004)

\bibitem{RARP2} Rodr\'{\i}guez, J. M., Alvarez, V., Romera, E., Pestana, D.:
Generalized weighted Sobolev spaces and
applications to Sobolev orthogonal polynomials II.
Approx. Theory and its Appl. 18:2, 1-32 (2002)

\bibitem{R1} Rodr\'{\i}guez, J. M.: Weierstrass' Theorem in
weighted Sobolev spaces. J. Approx. Theory 108, 119-160 (2001)

\bibitem{R2} Rodr\'{\i}guez, J. M.: The multiplication operator in
Sobolev spaces with respect to measures. J. Approx. Theory 109, 157-197 (2001)

\bibitem{R3} Rodr\'{\i}guez, J. M.: Approximation by polynomials
and smooth functions in Sobolev spaces with respect to measures.
J. Approx. Theory 120, 185-216 (2003)

\bibitem{R4} Rodr\'{\i}guez, J. M.:
A simple characterization of weighted Sobolev spaces with bounded multiplication operator.
To appear in J. Approx. Theory.

\bibitem{RY} Rodr\'{\i}guez, J. M., Yakubovich, D. V.:
A Kolmogorov-Szeg\"o-Krein type condition for weighted Sobolev spaces.
Indiana Univ. Math. J. 54, 575-598 (2005)



\end{thebibliography}
\end{document}